\documentclass{gtpart}     % Basic GT/GTM/AGT style
\gtart
\usepackage{latexsym}
\usepackage{amsfonts}
\usepackage{graphics}
\usepackage[all]{xy}
\usepackage{amscd}
\usepackage{multirow}
\title{Some new examples with almost positive curvature}

\author{Martin Kerin}
\givenname{Martin}
\surname{Kerin}
\address{Mathematisches Institut\\ WWU M\"unster\\ Einsteinstr. 62\\ 48149 M\"unster\\ Germany}
\email{m.kerin@math.uni-muenster.de}
\urladdr{}


\arxivreference{0810.4524}
\arxivpassword{da6j9}

\volumenumber{}
\issuenumber{}
\publicationyear{}
\papernumber{}
\startpage{}
\endpage{}
\doi{}
\MR{}
\Zbl{}
\received{}
\revised{}
\accepted{}
\published{}
\publishedonline{}
\proposed{}
\seconded{}
\corresponding{}
\editor{}
\version{}

\newtheorem{thm}{Theorem}[section]
\newtheorem*{thm*}{\bf Theorem}
\newtheorem{cor}[thm]{\bf Corollary}
\newtheorem*{cor*}{\bf Corollary}
\newtheorem{prop}[thm]{\bf Proposition}
\newtheorem*{prop*}{\bf Proposition}
\newtheorem{lem}[thm]{\bf Lemma}
\newtheorem*{lem*}{\bf Lemma}

\newtheorem*{defn*}{\bf Definition}

\newtheorem{main}{Theorem}                 %Theorems in introduction

\newtheorem*{qn*}{\bf Question}
\newtheorem*{conj*}{\bf Deformation Conjecture}

\newcommand{\pf}[1]{\noindent \emph{Proof of Theorem #1.}}

\theoremstyle{remark}
\newtheorem{rmk}[thm]{Remark}

\newtheorem{Rmk}{Remark}

\newtheorem{ex}[Rmk]{Example}

\def\3{{\ss}}
\def\2{\frac{1}{2}}
\def\4{\frac{1}{4}}
\def\8{\frac{1}{8}}
\def\x{\times}
\def\.{\cdot}

\def\ox{\otimes}
\def\bq{/\!\!/}
\def\<{\langle}
\def\>{\rangle}
\def\met{\<\; , \: \>}
\def\mmet{\<\!\<\; , \: \>\!\>}
\def\cmet{(\!(\; , \: )\!)}
\def\^{\wedge}
\def\vsig{\check \sigma}

\newcommand{\vphi}{\varphi}

\newcommand{\wt}{\widetilde}

\newcommand{\tr}{\ensuremath{\operatorname{tr}}}
\newcommand{\diag}{\ensuremath{\operatorname{diag}}}

\newcommand{\rank}{\ensuremath{\operatorname{rank}}}
\renewcommand{\Im}{\ensuremath{\operatorname{Im}}}
\newcommand{\Sec}{\ensuremath{\operatorname{sec}}}

\def\Hom{\mathop{\rm Hom}\nolimits}
\def\Im{\mathop{\rm Im}\nolimits}
\def\Re{\mathop{\rm Re\,}\nolimits}

\def\Ad{\mathop{\rm Ad}\nolimits}
\def\Span{\mathop{\rm Span}\nolimits}

\def\ceq{\mathrel{\mathop:}=}

\def\HH{\mathord{\mathbb H}}
\def\C{\mathord{\mathbb C}}
\def\R{\mathord{\mathbb R}}
\def\Z{\mathord{\mathbb Z}}
\def\NN{\mathord{\mathbb N}}
\def\Q{\mathord{\mathbb Q}}

\def\g{\mathfrak{g}}

\def\p{\mathfrak{p}}

\def\u{\mathfrak{u}}
\def\k{\mathfrak{k}}

\def\h{\mathfrak{h}}
\def\m{\mathfrak{m}}
\newcommand{\mf}{\mathfrak}

\def\L{\mathrm{L}}
\newcommand{\mcal}{\mathcal}

\newcommand{\CP}{\mathbf{CP}}
\newcommand{\Ca}{\mathbf{Ca}}
\newcommand{\rp}{\mathbf{RP}}

\newcommand{\sph}{\mathbf{S}}

\newcommand{\G}{\mathrm{G}}
\newcommand{\E}{\mathrm{E}}
\newcommand{\B}{\mathrm{B}}
\newcommand{\M}{\mathrm{M}}
\newcommand{\N}{\mathrm{N}}
\newcommand{\RR}{\mathrm{R}}
\newcommand{\T}{\mathrm{T}}
\newcommand{\GG}{\mathrm{GG}}
\newcommand{\K}{\mathrm{K}}
\def\P{\mathrm{P}}
\newcommand{\UU}{\mathrm{U}}
\def\H{\mathrm{H}}
\newcommand{\SU}{\operatorname{SU}}

\newcommand{\SO}{\operatorname{SO}}
\def\O{\operatorname{O}}
\newcommand{\Spin}{\operatorname{Spin}}
\newcommand{\U}{\operatorname{U}}
\newcommand{\pp}{\operatorname{p}}

\def\lra{\longrightarrow}
\def\lmt{\longmapsto}
\def\hra{\hookrightarrow}

\def\bsm{\begin{smallmatrix}}
\def\esm{\end{smallmatrix}}
\def\bpm{\begin{pmatrix}}
\def\epm{\end{pmatrix}}

\newcommand{\fourm}[4]{\mbox{$
\bpm
#1 &  &  & \\
   & #2 &  &  \\
   &  & #3 &  \\
  &  &  & #4 
\epm$}}

\def\beq{\begin{equation}}
\def\eeq{\end{equation}}
\def\Beq{\begin{equation*}}
\def\Eeq{\end{equation*}}
\def\beqs{\begin{eqnarray}}
\def\eeqs{\end{eqnarray}}
\def\Beqs*{\begin{eqnarray*}}
\def\Eeqs*{\end{eqnarray*}}

\begin{document}

\begin{abstract}    % type your abstract below
As a means to better understanding manifolds with positive curvature, there has been much recent interest in the study of non-negatively curved manifolds which contain points at which all $2$--planes have positive curvature.  We show that there are generalisations of the well-known Eschenburg spaces and quotients of $\sph^7 \x \sph^7$ which admit metrics with this property.
\end{abstract}

\maketitle

%%%%%%%%%%%%%%%%%%%%   Start of main body of article

It is an unfortunate fact that for a simply connected manifold which admits a metric of non-negative curvature there are no known obstructions to admitting positive curvature.  While there exist many examples of manifolds with non-negative curvature, the known examples with positive curvature are very sparse (see Ziller \cite{Zi} for a comprehensive survey of both situations).  Other than the rank-one symmetric spaces there are isolated examples in dimensions $6,7,12,13$ and $24$ due to Wallach \cite{Wa} and Berger \cite{Ber}, and two infinite families, one in dimension $7$ (Eschenburg spaces; see Aloff and Wallach \cite{AW}, Eschenburg \cite{E1}, \cite{E2}) and the other in dimension $13$ (Bazaikin spaces; see Bazaikin \cite{Ba}).  In recent developments, two distinct metrics with positive curvature on a particular cohomogeneity-one manifold have been proposed (Grove, Verdiani and Ziller \cite{GVZ}, Dearricott \cite{D}), while in \cite{PW2} Petersen and Wilhelm propose that the Gromoll--Meyer exotic $7$--sphere admits positive curvature, which would be the first exotic sphere known to exhibit this property.

In this paper we are interested in the study of manifolds which lie ``between'' those with non-negative and those with positive sectional curvature.  It is hoped that the study of such manifolds will yield a better understanding of the differences between these two classes.  Recall that a Riemannian manifold $(\M, \met)$ is said to have \emph{quasi-positive curvature} (resp. \emph{almost positive curvature}) if $(\M, \met)$ has non-negative sectional curvature and there is a point (resp. an open dense set of points) at which all $2$--planes have positive sectional curvature.

\begin{main}
  \label{thmA} \
  \begin{enumerate}
  \item \label{QPnew1}
    There exists a free circle action and a free $\sph^3$ action on $\sph^7 \x \sph^7$ such that each of the respective quotients $\M^{13} \ceq \sph^1 \backslash (\sph^7 \x \sph^7)$ and $\N^{11} \ceq \sph^3 \backslash (\sph^7 \x \sph^7)$ admits a metric with almost positive curvature.
  \item \label{QPnew2}
    If $\M^{13}$ and $\N^{11}$ are equipped with the metrics from \ref{QPnew1}, then there exist Riemannian submersions $\M^{13} \lra \CP^3$ and $\N^{11} \lra \sph^4$ such that in each case the fibre is $\sph^7$ and the bundle is non-trivial but shares the same cohomology ring as the corresponding product.
\end{enumerate}
\end{main}

We use the Pontrjagin class to distinguish $\M^{13}$ and $\N^{11}$ from the respective products in Theorem \ref{thmA}\ref{QPnew2}.  Moreover, in each case the induced metric on the base is positively curved.

It has been conjectured by F Wilhelm that, if $\M$ is a positively curved Riemannian manifold, then the dimension of the fibre of a Riemannian submersion $\M \lra \B$ must be less than the dimension of the base.  Theorem \ref{thmA} shows that this is false when the hypothesis is weakened to almost positive curvature.

Other than the Gromoll-Meyer exotic $7$--sphere (Gromoll and Meyer \cite{GM}, Wallach \cite{W}, Eschenburg and Kerin \cite{EK}, Petersen and Wilhelm \cite{PW2}), the only other previously known examples of manifolds with almost positive or quasi-positive curvature are given in Petersen and Wilhelm \cite{PW1}, Wilking \cite{Wi}, Tapp \cite{Ta1} and Kerin \cite{Ke2}.

While Wilking \cite{Wi} has shown that it is not possible in general to deform quasi-positive curvature to positive curvature, it is still unknown whether this can be achieved in the simply connected case or whether one can always deform quasi-positive curvature to almost positive curvature.

\begin{main}
  \label{thmB}
    Let $\L_{p,q} \subset \U(n+1) \x \U(n+1)$, $n \geq 2$, be defined by
            $$\qquad \L_{p,q} = \{(\diag(z^{p_1}, \dots, z^{p_{n+1}}), \diag(z^{q_1}, z^{q_2}, A)) \ | \ z \in \sph^1, A \in \U(n-1)\}$$
    where $p_1, \dots, p_{n+1}, q_1, q_2 \in \Z$.  When the two-sided action of $\L_{p,q}$ on $\U(n+1)$ is free, denote the quotient $\U(n+1) \bq \L_{p,q}$ by $\E^{4n-1}_{p,q}$.  All $\E^{4n-1}_{p,q}$ admit a metric with quasi-positive curvature.
\end{main}

Observe that to allow $p_i = p_j$ for all $1 \leq i < j \leq n+1$ in Theorem \ref{thmB} is equivalent to setting $p_i = 0$ for all $i$, since the centre of $\U(n+1)$ is given by multiples of the identity.  Wilking \cite{Wi} showed that these homogeneous spaces admit a metric with almost positive curvature whenever $q_1 q_2 < 0$, while Tapp \cite{Ta1} subsequently showed that $(q_1, q_2) \neq (0,0)$ is sufficient to guarantee the existence of a metric with quasi-positive curvature.

The biquotients $\E^{4n-1}_{p,q}$ in Theorem \ref{thmB} should be thought of as generalisations of the Eschenburg spaces, which arise when $n=2$.  In Eschenburg \cite{E1} it is shown that infinitely many Eschenburg spaces admit positive curvature, while in Kerin \cite{Ke2} it is shown that all Eschenburg spaces admit a metric with quasi-positive curvature.

The paper is organised as follows.  In Section \ref{Biqs} we review some notation and geometric techniques for biquotients.  In Section \ref{CayG2} we review some facts about the Cayley numbers and the exceptional Lie group ${\G}_2$.  In Section \ref{SOFree} we describe the manifolds $\M^{13}$ and $\N^{11}$ of Theorem \ref{thmA} as biquotients and as the total spaces of Riemannian submersions.  We prove the curvature statements of Theorem \ref{thmA}\ref{QPnew1} in Section \ref{QPMN}, while proof of the topological statements in \ref{QPnew2} is postponed until Section \ref{TopMN}.  Section \ref{GenEsch} is devoted to establishing Theorem \ref{thmB}.

%\ms
\emph{Acknowledgments:}  The majority of this work was completed as part of a Ph.D. thesis at the University of Pennsylvania under the supervision of Wolfgang Ziller.  His constant advice, support and friendship were inspiring.  Thanks also to Jost Eschenburg for useful comments.  Finally, I wish to thank Burkhard Wilking for suggesting some improvements.

\section{Biquotient actions and metrics}
\label{Biqs}

In his Habilitation, \cite{E1}, Eschenburg studied biquotients in great detail.  The following section provides a review of the material in \cite{E1} and establishes the basic language, notation and results which will be used throughout the remainder of the paper.

Let $\G$ be a compact Lie group, $\UU \subset \G \x \G$ a closed subgroup, and let $\UU$ act on $\G$ via
        $$(u_1, u_2) \star g = u_1 g u_2^{-1}, \ \  g \in \G, (u_1, u_2) \in \UU.$$
The action is free if and only if, for all non-trivial $(u_1, u_2) \in \UU$, $u_1$ is never conjugate to $u_2$ in $\G$.  The resulting manifold is called a {\it biquotient}.

Let $\K \subset \G$ be a closed subgroup, $\met$ be a left-invariant, right $\K$--invariant metric on $\G$, and $\UU \subset \G \x \K \subset \G \x \G$ act freely on $\G$ as above. Let $g \in \G$.  Define
        $$\UU^g_L \ceq \{(g u_1 g^{-1}, u_2) \ | \ (u_1, u_2) \in \UU \}.$$
Since $\UU$ acts freely on $\G$, so too does $\UU^g_L$, and $\G \bq \UU$ is isometric to $\G \bq \UU^g_L$.  This follows from the fact that left-translation $L_g \co \G \lra \G$ is an isometry which satisfies $g u_1 g^{-1}( L_g g') u_2^{-1} = L_g (u_1 g' u_2^{-1})$.  Therefore $L_g$ induces an isometry of the orbit spaces $\G\bq \UU$ and $\G\bq \UU^g_L$.

Consider a Riemannian submersion $\pi\co \M^n \lra \N^{n-k}$.  By O'Neill's formula for Riemannian submersions, $\pi$ is curvature non-decreasing.  Therefore $\Sec_\M \geq 0$ implies $\Sec_\N \geq 0$, and zero--curvature planes on $\N$ lift to horizontal zero--curvature planes on $\M$.  Because of the Lie bracket term in the O'Neill formula the converse is not true in general, namely horizontal zero--curvature planes in $\M$ cannot be expected to project to zero--curvature planes on $\N$.

Let $\K \subset \G$ be Lie groups, $\k \subset \g $ the corresponding Lie algebras, and $\met$ a non-negatively curved left-invariant metric on $\G$ which is right-invariant under $\K$.  We can write $\g = \k \oplus \p$ with respect to $\met$.  Given $X \in \g$ we will always use $X_\k$ and $X_\p$ to denote the $\k$ and $\p$ components of $X$ respectively.
\begin{gather*}\tag*{\text{Recall that}}
        \G \cong (\G \x \K)/ \Delta \K
\end{gather*}
via $(g,k) \lmt g k^{-1}$, where $\Delta \K$ is the free diagonal action of $\K$ on the right of $\G \x \K$.  Notice that the restriction of $\met$ to $\K$ is bi-invariant.  Equip $\G \x \K$ with the metric $\mmet \ceq \met \oplus t \met |_\k$, $t > 0$.  Thus we may define a new left-invariant, right $\K$--invariant metric $\met_1$ (with $\Sec \geq 0$) on $\G$ via the Riemannian submersion
        \begin{gather}
        \begin{aligned}
        \label{newmet}
            (\G \x \K, \mmet) &\lra (\G, \met_1)\\
            (g,k) &\lmt g k^{-1}
        \end{aligned}\\
        \tag*{\text{with}}
            \< X, Y \>_1 = \< X, \Phi(Y) \>
        \end{gather}%\eeq
%        \begin{align}
%        \nonumber
%            (\G \x \K, \mmet) &\lra (\G, \met_1)\\
%        \label{newmet}
%            (g,k) &\lmt g k^{-1}
%        \end{align}
%        \begin{gather*}\tag*{\text{with}}
%        %\beq
%            \< X, Y \>_1 = \< X, \Phi(Y) \>
%        \end{gather*}%\eeq
where $\Phi(Y) = Y_\p + \lambda Y_\k$, $\lambda = \frac{t}{t+1} \in (0,1)$.  Furthermore, it is clear that the metric tensor $\Phi$ is invertible with inverse described by $\Phi^{-1}(Y) = Y_\p + \frac{1}{\lambda} Y_\k$.

Suppose $\sigma = \Span \{\Phi^{-1}(X), \Phi^{-1}(Y)\} \subset \g$ is a zero--curvature plane with respect to the metric $\met_1$, i.e. $\Sec_1(\sigma) = 0$.  By the O'Neill formula $\sigma$ must therefore lift to a horizontal zero--curvature plane $\wt \sigma \subset \g \oplus \k$ with respect to $\mmet$.  It is easy to check that the horizontal lift of a vector $\Phi^{-1}(X) \in \g$ to $\g \oplus \k$ is given by $\left(X, -\frac{1}{t} X_\k\right)$.  Then clearly
        $$\wt \sigma = \Span\left\{\left(X, -\frac{1}{t} X_\k\right), \left(Y, -\frac{1}{t} Y_\k\right) \right\}.$$
But, since $\mmet$ is a non-negatively curved product metric, it follows immediately by considering the unnormalised curvature that $\wt \sigma$ has zero--curvature if and only if $[X_\k, Y_\k] = 0$ and the plane $\Span\{X,Y\} \subset \g$ has zero--curvature with respect to the original metric $\met$, i.e. $\Sec(X,Y) = 0$.

From Tapp \cite{Ta2}, which generalizes similar results in Eschenburg \cite{E1} and Wilking \cite{Wi}, we know that if $\met$ is induced by a Riemannian submersion to $\G$ from a bi-invariant metric on some Lie group $\L$, then in fact $\Sec_1(\sigma) = 0$ if and only if $\Sec(\wt \sigma) = 0$ with respect to $\mmet$, i.e. if and only if $\Sec(X,Y) = 0$ and $[X_\k, Y_\k] = 0$.  We will always be in this situation as throughout the paper we will use only the metrics described in Examples \ref{Eschlem1} and \ref{Eschlem2} below.

\begin{ex}
\label{Eschlem1}
Suppose that $(\G,\K)$ is a symmetric pair and that the initial metric $\met = \met_0$ is a bi-invariant metric on $\G$.  As in (\ref{newmet}), equip $\G$ with a new metric
        \beq
        \label{met1}
            \< X, Y \>_1 = \< X, \Phi_1(Y) \>_0
        \eeq
where $\g = \k \oplus \p$ with respect to $\met_0$ and $\Phi_1(Y) = Y_\p + \lambda_1 Y_\k$.  Then $\sigma = \Span \{\Phi_1^{-1}(X),\Phi_1^{-1}(Y)\} \subset \g$ has zero--curvature with respect to $\met_1$, i.e. $\Sec_1(\sigma) = 0$, if and only if
        \beq
        \label{met1conds}
            0 = [X,Y] = [X_\k,Y_\k] = [X_\p,Y_\p].
        \eeq
The proof of this follows immediately from our previous discussion together with the fact that $[\p, \p] \subset \k$ whenever $(\G,\K)$ is a symmetric pair.
\end{ex}

\begin{ex}
\label{Eschlem2}
Let $\G \supset \K \supset \H$ be a chain of subgroups and suppose that both $(\G,\K)$ and $(\K,\H)$ are symmetric pairs. Let $\g = \k \oplus \p$ and $\k = \h \oplus \m$ be the corresponding orthogonal decompositions with respect to the bi-invariant metric $\met_0$ on $\G$.  Start with the metric $\met = \met_1$ defined by Example \ref{Eschlem1}.  Now define the metric $\met_2$ on $\G$ as in (\ref{newmet}), where $\K$ is replaced by $\H$, $s>0$ takes the role of $t$, and $\Psi$ replaces $\Phi$:
        \begin{align}
        \label{met2met1}
            \<X,Y\>_2 &= \<X,\Psi(Y)\>_1 \\
        \label{met2met0}
                &= \<X,\Phi_2(Y)\>_0
        \end{align}
with $\Phi_2(Y) = Y_\p + \lambda_1 Y_\m + \lambda_1 \lambda_2 Y_\h$, $\lambda_2 = \frac{s}{s+1}$, and $\Psi(Y) = \Phi_1^{-1} \Phi_2(Y) = Y_\p + Y_\m + \lambda_2 Y_\h$.

Let $\sigma = \Span \{\Psi^{-1}(X),\Psi^{-1}(Y)\} \subset \g$.  Then, by our discussion prior to Example \ref{Eschlem1}, $\Sec_2(\sigma) = 0$ if and only if
$\Sec_1(X,Y) = 0$ and $[X_\h, Y_\h] = 0$.  By again considering horizontal lifts it is not difficult to check that $\Sec_1(X,Y) = 0$ if and only if conditions (\ref{met1conds}) hold as for $\Sec_1(\Phi_1^{-1} (X),\Phi_1^{-1} (Y)) = 0$.  Hence $\Sec_2(\sigma) = 0$ if and only if
        \beq
        \label{met2conds}
            0 = [X, Y] = [X_\k, Y_\k] = [X_\p,Y_\p] = [X_\m,Y_\m] = [X_\h,Y_\h]
        \eeq
where we have used the fact that $[\m,\m] \subset \h$ since $(\K, \H)$ is a symmetric pair.
\end{ex}

Now that we have described how to induce new metrics on $\G$ from old ones and derived zero--curvature conditions for these metrics, we proceed to consider biquotients $\G \bq \UU$. Let $\Delta \G = \{(g,g) \ | \ g \in \G\}$.  Then, if the two-sided action of $\UU$ on $\G$ is free, $\Delta \G \x \UU$ acts freely on $\G \x \G$ via
        \beq
        \label{GUaction}
            ((g, g), (u_1, u_2)) \star (g_1, g_2) = (g g_1 u_1^{-1}, g g_2 u_2^{-1})
        \eeq
with $(g, g) \in \Delta \G$, $(u_1, u_2) \in \UU$, $(g_1, g_2) \in \G \x \G$, and there is a canonical diffeomorphism
        \beq
        \label{doubling}
            \Delta \G \backslash (\G \x \G)/ \UU \cong \G \bq \UU
        \eeq
induced by the map
        \begin{align*}
            \G \x \G &\lra \G\\
            (g_1,g_2) &\lmt g_1^{-1} g_2.
        \end{align*}
Let $\K_1$ and $\K_2$ be arbitrary subgroups of $\G$.  We define left-invariant metrics, $\met_{\K_1}$ and $\met_{\K_2}$, on $\G$ as in (\ref{newmet}).  Equip $\G \x \G$ with a left-invariant, right $(\K_1 \x \K_2)$--invariant product metric $\cmet = \met_{\K_1} \oplus \met_{\K_2}$.  If $\UU \subset \K_1 \x \K_2$ then the $\Delta \G \x \UU$ action is by isometries and $\cmet$ induces a metric on $\G \bq \UU$.  Our goal is to determine when a plane tangent to $\G \bq \UU$ has zero--curvature with respect to this induced metric.

By (\ref{doubling}) and our choice of metric, the quotient map $(\G \x \G, \cmet) \lra \G \bq \UU$ is a Riemannian submersion.  O'Neill's formula implies that a zero--curvature plane tangent to $\G \bq \UU$ must lift to a horizontal zero--curvature plane with respect to $\cmet$.  As in the case of metrics on $\G$, if $\cmet$ is induced from a bi-invariant metric on some Lie group $\L$, then Tapp \cite{Ta2} implies that horizontal zero--curvature planes with respect to $\cmet$ must project to zero--curvature planes in $\G \bq \UU$.  For our purposes this will always be true since we will consider only metrics as in Examples \ref{Eschlem1} and \ref{Eschlem2}.

We must determine what it means for a plane to be horizontal with respect to $\cmet$ and the $\Delta \G \x \UU$ action.  Since each $\Delta \G \x \UU$ orbit passes through some point of the form $(g, e) \in \G \x \G$, where $e$ is the identity element of $\G$, we may restrict our attention to such points.

Recall that $\cmet$ is left-invariant.  Therefore, letting $\u$ denote the Lie algebra of $\UU$, the vertical subspace at $(g,e) \in \G \x \G$ is given by
        $$\mcal{V}_g = \left\{ \left(\Ad_{g^{-1}} X - Y_1, X - Y_2 \right) \ | \ X \in \g, (Y_1, Y_2) \in \mf{u} \right\}$$
after left-translation to $(e,e) \in \G \x \G$.  Note that this is independent of the choice of left-invariant metric on $\G \x \G$.

Thus, with respect to $\cmet$, the horizontal subspace at $(g, e)$ is
        \beq
        \label{horizsp}
            \mcal{H}_g = \left\{ \left(\Omega_1^{-1}(- \Ad_{g^{-1}} X), \Omega_2^{-1}(X) \right) \ | \ \<X, \Ad_g Y_1 - Y_2\>_0 = 0 \ \forall \ (Y_1, Y_2) \in \mf{u} \right\}
        \eeq
where $\Omega_1$ and $\Omega_2$ are the metric tensors relating the left-invariant metrics $\met_{\K_1}$ and $\met_{\K_2}$ respectively to a fixed bi-invariant metric $\met_0$ on $\G$, i.e. $\<X, Y\>_{\K_i} = \<X, \Omega_i(Y)\>_0$, $i=1,2$.  We recall that the metric tensors in Examples \ref{Eschlem1} and \ref{Eschlem2} are given by $\Phi_1$ and $\Phi_2$ respectively, as shown in (\ref{met1}) and (\ref{met2met0}).

In particular, (\ref{horizsp}) shows that a horizontal $2$--plane $\sigma$ in $(\G \x \G, \cmet)$ must project to a $2$--plane on each factor, denoted by $\vsig_1$ and $\vsig_2$ respectively. Moreover, since $\cmet$ is a product metric, $\Sec(\sigma) = 0$ if and only if $\Sec_1 (\vsig_i) = 0$, $i = 1,2$.  Thus, for product metrics involving the metrics described by Examples \ref{Eschlem1} and \ref{Eschlem2}, we may apply conditions (\ref{met1conds}) and (\ref{met2conds}) respectively in order to determine when a horizontal plane $\sigma$ has zero--curvature.

\section{The Cayley numbers, $\G_{2}$ and its Lie algebra}
\label{CayG2}

We recall without proof some well known facts about Cayley numbers, the Lie group $\G_2$ and its Lie algebra.  More details may be found in Gluck, Warner and Ziller \cite{GWZ} and Murakami \cite{M}.

We may write the Cayley numbers as $\Ca = \HH + \HH \ell$, where $\HH$ is the algebra of quaternions.  Thus we have a natural orthonormal basis
                 $$\{e_0 = 1, e_1 = i, e_2 = j, e_3 = k, e_4 = \ell, e_5 = i \ell, e_6 = j \ell, e_7 = k \ell\}$$
for $\Ca$.  Note that this description of $\Ca$ differs slightly from that given in Murakami \cite{M}, and accounts for the difference which occurs in the description of the Lie algebra $\g_2$ in Theorem \ref{g2eqns}.  Multiplication in $\Ca$ is non-associative and defined via
        \beq
        \label{multiplicationCa}
            (a+b\ell)(c+d\ell) = (ac - \bar d b)+(da + b \bar c)\ell, \ \ a,b,c,d \in \HH.
        \eeq
Hence we have the following multiplication table, where the order of multiplication is given by (row)$\ast$(column):
        \begin{table}[!hbtp]
        \begin{center}
        \begin{tabular}{|c||c|c|c|c|c|c|c|} \cline{2-8}
            \multicolumn{1}{c||}{} & $e_1 = i$ & $e_2 = j$ & $e_3 = k$ & $e_4 = \ell$ & $e_5 = i \ell$ & $e_6 = j \ell$ & $e_7 = k \ell$ \\ \hline \hline
            $e_1 = i$ & $-1$ & $k$ & $-j$ & $i \ell$ & $-\ell$ & $-k \ell$ & $j \ell$  \\ \hline
            $e_2 = j$ & $-k$ & $-1$ & $i$ & $j \ell$ & $k \ell$ & $-\ell$ & $-i \ell$    \\ \hline
            $e_3 = k$ & $j$ & $-i$ & $-1$ & $k \ell$ & $-j \ell$ & $i \ell$ & $-\ell$  \\ \hline
            $e_4 = \ell$ & $-i \ell$ & $-j \ell$ & $-k \ell$ & $-1$ & $i$ & $j$ & $k$ \\ \hline
            $e_5 = i \ell$ & $\ell$ & $-k \ell$ & $j \ell$ & $-i$ & $-1$ & $-k$ & $j$ \\ \hline
            $e_6 = j \ell$ & $k \ell$ & $\ell$ & $-i \ell$ & $-j$ & $k$ & $-1$ & $-i$ \\ \hline
            $e_7 = k \ell$ & $-j \ell$ & $i \ell$ & $\ell$ & $-k$ & $-j$ & $i$ & $-1$ \\ \hline

        \end{tabular}
        \end{center}
        \caption{Multiplication table for $\Ca$} \label{multiplicationtable}
        \end{table}

Recall that the Lie group $\G_2$ is the automorphism group of $\Ca \cong \R^8$.  In fact $\G_2$ is a connected subgroup of $\SO(7) \subset \SO(8)$, where $\SO(8)$ acts on $\Ca \cong \R^8$ by orthogonal transformations and $\SO(7)$ is that subgroup consisting of elements which leave $e_0 = 1$ fixed.  $\SO(8)$ also contains two copies of $\Spin(7)$ which are not conjugate in $\SO(8)$, and $\G_2$ is the intersection of these two subgroups.

As our eventual goal is to prove Theorem \ref{thmA}\ref{QPnew1} and \ref{QPnew2}, it is useful to recall the fact that $\G_2$ appears in the descriptions of some interesting homogeneous spaces.  The following statements are well-known and follow from applications of the triality principle for $\SO(8)$.  More details may be found in, for example, Adams \cite[Theorem 5.5]{Ad}, Murakami \cite{M}, Jacobson \cite[page 79]{J}.

\begin{thm} \
\label{G2HomogSp}
\begin{enumerate}
    \item
    \label{G2sphere}
        $\Spin(7)/\G_2 = \sph^7$, which inherits constant positive curvature from the bi-invariant metric on $\Spin(7)$.  Moreover, the isotropy representation is transitive on the collection of pairs of orthogonal unit tangent vectors.
    \item
    \label{S7xS7}
        $\Spin(8)/\G_2 = \sph^7 \x \sph^7$ and $\SO(8)/\G_2 = (\sph^7 \x \sph^7)/ \Z_2$, where $\Z_2 = \{\pm ({\rm id}, {\rm id}) \}$.
    \item
    \label{S6}
        $\G_2/\SU(3) = \sph^6$.
\end{enumerate}
\end{thm}

We now turn our attention to the Lie algebra of $\G_2$.  The proof of the following theorem follows exactly as in Murakami \cite{M} except that we use the basis and multiplication conventions for $\Ca$ as in Table \ref{multiplicationtable}.  Recall that $\mf{so}(n) = \{ A \in M_n (\R) \ | \ A^t = -A\}$.

\begin{thm}
\label{g2eqns}
The Lie algebra of $\G_2$, denoted by $\g_2$, consists of matrices $A = (a_{ij}) \in \mf{so}(7)$ which satisfy $a_{ij} + a_{ji} = 0$ and the following equations:
    \begin{align*}
        a_{23} + a_{45} +a_{76} &= 0\\
        a_{12} + a_{47} + a_{65} &= 0\\
        a_{13} + a_{64} + a_{75} &= 0\\
        a_{14} + a_{72} + a_{36} &= 0\\
        a_{15} + a_{26} + a_{37} &= 0\\
        a_{16} + a_{52} + a_{43} &= 0\\
        a_{17} + a_{24} + a_{53} &= 0
    \end{align*}
Hence $\g_2 \subset \mf{so}(7)$ is $14$--dimensional and consists of matrices of the form:
    \beq
    \label{g2}
        \bpm
            0 & x_1 + x_2  & y_1 + y_2 & x_3 + x_4 & y_3 + y_4  & x_5 + x_6 \! & y_5 + y_6 \\
            -(x_1 + x_2)  & 0 & \alpha_1 & -y_5 & x_5  & -y_3  & x_3 \\
            -(y_1 + y_2)  & -\alpha_1 & 0 & x_6 & y_6  & -x_4  & -y_4 \\
            -(x_3 + x_4)  & y_5 & -x_6 & 0 & \alpha_2  & y_1  & -x_1 \\
            -(y_3 + y_4)  & -x_5 & -y_6 & -\alpha_2 & 0  & x_2  & y_2 \\
            -(x_5 + x_6)  & y_3 & x_4 & -y_1 & -x_2  & 0  & \alpha_1 + \alpha_2 \\
            -(y_5 + y_6)  & -x_3 & y_4 & x_1 & -y_2  & -(\alpha_1 + \alpha_2)  & 0
        \epm
    \eeq
\end{thm}
Recall that $\G_2$ is a rank $2$ Lie group.  Therefore an examination of the elements (\ref{g2}) of $\g_2$ reveals that the maximal torus of $\G_2$ is given by:
        \beq
        \label{maxtorus}
            \T^2 = \left\{    \fourm{1}{R(\theta)}{R(\varphi)}{R(\theta + \varphi)}
                     \ \Bigg| \ R(\theta) = \bpm \cos \theta & - \sin \theta \\ \sin \theta & \cos \theta \epm \right\}
        \eeq

\section{Free isometric actions on $\mathrm{SO}(8)$}
\label{SOFree}

Consider the rank one symmetric pair $(\G, \K) = (\SO(8), \SO(7))$ where
        \begin{align*}
            \SO(7) &\hra \SO(8)\\
            A &\lmt \bpm 1 & \\ & A \epm
        \end{align*}
with Lie algebras $\g, \k$ respectively.  Let $\< X, Y \>_0 = - \tr (XY)$ be a bi-invariant metric on $\G$.  With respect to $\met_0$ we thus have $\g = \p \oplus \k$.  As in (\ref{met1}) we define a left-invariant, right $\K$--invariant metric $\met_1$ on $\G$ by
        \beq
            \< X,Y\>_1 = \<X, \Phi_1(Y)\>_0
        \eeq
where $\Phi_1(Y) = Y_\p + \lambda_1 Y_\k$, $\lambda_1 \in (0,1)$.  Recall that from Example \ref{Eschlem1} we know that a plane
        $$\sigma = \Span\{\Phi_1^{-1}(X), \Phi_1^{-1}(Y)\} \subset \g$$
has zero--curvature with respect to $\met_1$ if and only if
        \beq
        \label{EschLemConds}
            0 = [X, Y] = [X_\p, Y_\p] = [X_\k, Y_\k].
        \eeq
We now equip $\G \x \G$ with the product metric $\met_1 \oplus \met_0$.

Consider an isometric action of $\UU \ceq \sph^1 \x \G_2 \subset \K \x \G$ on $\SO(8)$ defined by
        \beq
        \label{actions}
            A \lmt \wt R (\theta) \. A \. g^{-1}
        \eeq
where $A \in \SO(8)$, $g \in \G_2$, and
       \beq
        \label{circle}
            \wt R (\theta) = \fourm{I_{2 \x 2}}{R(p_1 \theta)}{R(p_2 \theta)}{R(p_3 \theta)}, \ \
            R(\theta) = \bpm \cos \theta & - \sin \theta \\ \sin \theta & \cos \theta\epm.
        \eeq
From (\ref{doubling}) we know that $\Delta \G \backslash \G \x \G / \UU \cong \G \bq \UU$ whenever the biquotient action of $\UU$ on $\G$ is free.

\begin{lem}
\label{Freeactions}
    $\Delta \G \x \UU$ acts freely and isometrically on $(\G \x \G, \met_1 \oplus \met_0)$ if and only if $(p_1, p_2, p_3)$ is equal to $(0,0,1)$ (up to sign and permutations of the $p_i$).
\end{lem}

\begin{proof}
Recall that conjugation of either factor of $\UU$ by elements of $\G$ is a diffeomorphism, and that a biquotient action is free if and only if non-trivial elements in each factor are never conjugate to one another in $\G$.  Thus we need only show that non-trivial elements of $\sph^1$ and $\T^2$ are never conjugate in $\G$ if and only if $(p_1, p_2, p_3) = (0,0,1)$ up to sign and permutations of the $p_i$, where $\T^2$ is the maximal torus of $\G_2$ described in (\ref{maxtorus}).  This amounts to investigating when the sets of $2 \x 2$ blocks on each side are equal up to conjugation by an element of the Weyl group of $\SO(8)$.  We recall that the Weyl group of $\SO(2n)$ acts via permutations of the $2 \x 2$ blocks and changing an even number of signs, where by a change of sign we mean $R(\theta) \lmt R(- \theta)$.  A simple calculation then yields the result.
\end{proof}

Note that there are many other free $\sph^1 \x \G_2$ actions on $\G$.  For example, there is a free $\sph^1$ action on the left of $\G / \G_2$ by matrices of the form
        \beq
        \label{otheractions}
            \fourm{R(\theta)}{R(\theta)}{R(\theta)}{R(k \theta)}
        \eeq
where $(k, 3) = 1$.  However, it is clear that only the action in Lemma \ref{Freeactions} is isometric with respect to the metric $\met_1 \oplus \met_0$ on $\G \x \G$.

It follows immediately from the long exact homotopy sequence for fibrations that a biquotient $\Spin(8) \bq (\sph^1 \x \G_2) = \sph^1 \backslash (\sph^7 \x \sph^7)$ must be simply connected.  By the lifting criterion for covering spaces the action by $\UU$ on $\SO(8)$ described above lifts to some action by $\sph^1 \x \G_2$ on $\Spin(8)$. Therefore, together with Theorem \ref{G2HomogSp}, one might expect that the resulting simply connected biquotient $\Spin(8) \bq (\sph^1 \x \G_2) = \sph^1 \backslash (\sph^7 \x \sph^7)$ is a non-trivial finite cover of $\SO(8) \bq (\sph^1 \x \G_2)$.  In fact the lemma below will demonstrate that this covering map is a diffeomorphism.

\begin{lem}
\label{lucky}
$\M^{13} \ceq \SO(8) \bq (\sph^1 \x \G_2)$ is simply connected and hence a quotient of $\sph^7 \x \sph^7$ by an $\sph^1$ action.
\end{lem}

\begin{proof}
Consider a general embedding
        \begin{align*}
            \sph^1_q &\hra \SO(8)\\
            R(\theta) &\lmt \fourm{R(q_1 \theta)}{R(q_2 \theta)}{R(q_3 \theta)}{R(q_4 \theta)}
        \end{align*}
where $q = (q_1, q_2, q_3, q_4) \in \Z^4$, where $R(u) \in \SO(2)$.  The long exact homotopy sequence for the fibration $\sph^1_q \x \G_2 \lra \SO(8) \lra \SO(8) \bq \sph^1_q \x \G_2$ yields
        $$\dots \lra \pi_1 (\sph_q^1 \x \G_2) = \Z \lra \pi_1 (\SO(8)) = \Z_2 \lra \pi_1 (\SO(8) \bq \sph_q^1 \x \G_2) \lra 0.$$
Thus to obtain the desired result we need only show that the map $\Z \lra \Z_2$ is surjective.

Recall that the homomorphism $\iota_* \co \pi_1 (\sph_q^1) \lra \pi_1 (\SO(n))$ is determined by the weights $q = (q_1, \dots, q_m)$, $m = \lfloor \frac{n}{2} \rfloor$, of the embedding, namely $\iota_*(1) = \sum q_i \mod 2$.  Therefore $\iota_*$ is onto exactly when $\sum q_i$ is odd.  In our case we have $q = (0,0,0,1)$, and so $\iota_*$ is a surjection.
\end{proof}

Notice that the action of $\UU$ on $\SO(8)$ given in Lemma \ref{Freeactions} may be enlarged to an isometric action by $\SO(3) \x \G_2$, and the resulting biquotient we call $\N^{11}$.  Now recall that for all $n$ we have a $2$-fold cover $\Spin(n) \lra \SO(n)$ with $\pi_1 (\Spin(n)) = 0$ and $\pi_1 (\SO(n)) = \Z_2$.  Thus, by the lifting criterion for covering spaces, the inclusion $\SO(3) \hra \SO(8)$ must lift to $\Spin(3) = \sph^3 \hra \Spin(8)$.  As in the case of $\UU = \sph^1 \x \G_2$ above we show that $\N^{11} = \SO(8) \bq (\SO(3) \x \G_2)$ is simply connected and hence diffeomorphic to $\Spin(8) \bq \sph^3 \x \G_2 = \sph^3 \backslash (\sph^7 \x \sph^7)$.

\begin{lem}
\label{lucky2}
$\N^{11} = \SO(8) \bq (\SO(3) \x \G_2)$ is simply connected and hence a quotient of $\sph^7 \x \sph^7$ by an $\sph^3$ action.
\end{lem}

\begin{proof}
Consider the chain of embeddings $j \circ i \co \sph^1 = \SO(2) \hra \SO(3) \hra \SO(8)$ given by enlarging $\sph^1$ above to an $\SO(3)$ in $\SO(8)$.  We thus have an induced homomorphism on fundamental groups $(j \circ i)_* = j_* \circ i_* \co \Z \lra \Z_2 \lra \Z_2$.  But $i_*$ and $(j \circ i)_*$ are simply the homomorphism $\iota_*$ from Lemma \ref{lucky}.  Hence $i_*(1) = 1 \mod 2$ and $(j \circ i)_*(1) = 1 \mod 2$.  This implies $j_*(1) = 1 \mod 2$ and therefore $j_*$ is a surjection.  An examination of the long exact homotopy sequence of the fibration $\SO(3) \x \G_2 \lra \SO(8) \lra \N^{11}$ yields the result.
\end{proof}

Recall that $\G_2 \subset \Spin(7)' \subset \SO(8)$, where the second inclusion is via the spin embedding.  By our choice of metric on $\G \x \G$ we may therefore enlarge the actions of $\sph^1 \x \G_2$ and $\SO(3) \x \G_2$ to isometric actions by $\sph^1 \x \Spin(7)'$ and $\SO(3) \x \Spin(7)'$ respectively.  It was shown in Eschenburg \cite[Theorem 75, Table 101]{E1} that these actions are free and the resulting biquotients are $\CP^3$ and $\sph^4$ respectively.  It follows immediately that there are fibre bundles
\begin{align*}
&\sph^7 = \Spin(7)'/\G_2 \lra \M^{13} \lra \CP^3 \\
&\sph^7 = \Spin(7)'/\G_2 \lra \N^{11} \lra \sph^4
\end{align*}
for which the projections are Riemannian submersions.

\section{Almost positive curvature on $\M^{13}$ and $\N^{11}$}
\label{QPMN}

We are now in a position to prove the curvature statements Theorem \ref{thmA}\ref{QPnew1}.  We will concentrate on the circle quotient of $\sph^7 \x \sph^7$, namely
        $$\M^{13} = \SO(8) \bq (\sph^1 \x \G_2) = \G \bq \UU,$$
since the other case follows trivially.

Consider the inclusions $\G = \SO(8) \supset \K = \SO(7) \supset \G_2$.  With respect to the bi-invariant metric $\< X, Y \>_0 = - \tr (XY)$ on $\G$ we have
        $$\g = \p \oplus \k \ \ \ {\rm and} \ \ \ \k = \m \oplus \g_2$$
where
        \beq
        \label{p}
            \p = \left\{\left(\begin{array}{c|ccc}
                                 0 & & -w^t & \\
                            \hline & & & \\
                                 w & & 0 & \\
                                   & & & \end{array} \right) \ \Big| \ w \in \R^7 \right\}
        \eeq
and, by (\ref{g2}):
        \beq
        \label{m}
            \m = \left\{\left(\begin{array}{c|ccccccc}
                            0 & 0 & 0 & 0 & 0 & 0 & 0 & 0\\
                     \hline 0 & 0 & v_1 & v_2 & v_3 & v_4 & v_5 & v_6 \\
                            0 & -v_1 & 0 & v_7 & v_6 & -v_5 & v_4 & -v_3\\
                            0 & -v_2 & -v_7 & 0 & -v_5 & -v_6 & v_3 & v_4\\
                            0 & -v_3 & -v_6 & v_5 & 0 & v_7 & -v_2 & v_1\\
                            0 & -v_4 & v_5 & v_6 & -v_7 & 0 & -v_1 & -v_2\\
                            0 & -v_5 & -v_4 & -v_3 & v_2 & v_1 & 0 & -v_7\\
                            0 & -v_6 & v_3 & -v_4 & -v_1 & v_2 & v_7 & 0
                        \end{array}\right) \ \Bigg| \ v_i \in \R \right\}
        \eeq

\begin{lem}
\label{rklem}
Let $W_i \in \p$, $i = 1,2$, and $V_j \in \m$, $j = 1,2$, be orthonormal vectors with respect to $\met_0$.  Then $\rank(W_i) = 2$, $\rank(V_j) = 6$, $\rank([W_1, W_2]) = 2$ and $\rank([V_1, V_2]) = 6$.
\end{lem}

\begin{proof}
Suppose $P$, $Q$ and $Z$ are real $n \x n$ matrices.  It is well-known (and not difficult to prove using the rank-nullity theorem) that:
\begin{enumerate}
\item If $\rank(P) = n$, then $\rank(PZ) = \rank(Z)$
\item If $\rank(Q) = n$, then $\rank(ZQ) = \rank(Z)$.
\end{enumerate}
In particular, if $P$ and $Q$ are invertible matrices then, for $Z_0 = PZQ$, it follows that  $\rank(Z) = \rank(Z_0)$.

Choose arbitrary pairs of orthonormal vectors in $\p$ and $\m$ and show that the conclusions of the lemma hold.  The result  now follows from the facts that the isotropy representation of $\G_2 \subset \SO(8)$ splits irreducibly as $\p \oplus \m$ and, moreover, is transitive on the sets of pairs of orthonormal vectors in $\p$ and $\m$ respectively, i.e. $\G_2$ acts transitively on the unit tangent bundle to $\sph^6$.  (Recall Theorem \ref{G2HomogSp}\ref{G2sphere}.)
\end{proof}

Let $\UU = \sph^1 \x \G_2 \subset \K \x \G$ be as in Lemma \ref{Freeactions}.  Thus, equipping $\G \x \G$ with the product metric $ \met_1 \oplus \met_0$ as before, we may induce a metric on $\G \bq \UU$ via the diffeomorphism
        $$\Delta \G \backslash \G \x \G / \UU \lra \G \bq \UU.$$
As discussed in Section \ref{Biqs}, we may restrict our attention to points of the form $(A,I) \in \G \x \G$.  Let $E_{78} \in \g$ denote the vector spanning the Lie algebra of the $\sph^1$ factor of $\UU$.  By (\ref{horizsp}) the horizontal subspace at $(A,I)$ with respect to $ \met_1 \oplus \met_0$ is given by
        $$
\mcal{H}_A = \{(- \Phi_1^{-1} (\Ad_{A^t} W), W) \ | \ W_{\g_2} = 0, \<W, \Ad_A E_{78} \>_0 = 0\}.
        $$
Suppose that
$$
\sigma = \Span \{(- \Phi_1^{-1} (\Ad_{A^t} X), X), (- \Phi_1^{-1} (\Ad_{A^t} Y), Y)\} \subset \g \oplus \g
$$
is a horizontal zero--curvature plane at $(A, I) \in \G \x \G$.  Since we have equipped $\G \x \G$ with the product metric $ \met_1 \oplus \met_0$, our discussion in Section \ref{Biqs} shows that $\sigma$ must project to zero--curvature planes $\vsig_i$, $i = 1,2$, on each factor.  Considering $\vsig_2$ we find:

\begin{lem}
\label{XY}
     Suppose $X, Y \in \g$ are linearly independent vectors such that $X_{\g_2} = Y_{\g_2} = 0$ and $[X,Y] = 0$.  Then it may be assumed without loss of generality that $X \in \p$ and $Y \in \m$.
\end{lem}

\begin{proof}
Since $[X,Y] = 0$ and $(\G,\K)$ is a symmetric pair, it follows that $[X_\p, Y_\p] = - [X_\m, Y_\m]$.  Suppose that $[X_\m, Y_\m] \neq 0$.  Then, by Lemma \ref{rklem}, $[X_\m, Y_\m]$ is a matrix of rank $6$.  But $[X_\p, Y_\p]$ has rank no more than $2$, so we have a contradiction.  Therefore $[X_\p, Y_\p] = [X_\m, Y_\m] = 0$.

$[X_\p, Y_\p] = 0$ implies that, since $(\G,\K)$ is a rank one symmetric pair, we may assume  $Y_\p = 0$ without loss of generality.  Hence $X \in \p \oplus \m$ and $Y \in \m$.  On the other hand, by Theorem \ref{G2HomogSp}\ref{G2sphere} there are no independent commuting vectors in $\m$.  Then, without loss of generality, $X \in \p$ and $Y \in \m$.
\end{proof}

The horizontal zero--curvature plane $\sigma$ is thus determined by
        \beq \label{XYvec}
            X = \left(\begin{array}{c|ccc}
                                 0 & & -w^t & \\
                            \hline & & & \\
                                 w & & 0 & \\
                                   & & & \end{array} \right) \in \p, \ \ \
            Y = \left(\begin{array}{c|ccc}
                                 0 & & 0 & \\
                            \hline & & & \\
                                 0 & & (v_{ij}) & \\
                                   & & & \end{array} \right) \in \m
        \eeq
with $[X, Y] = 0$, where $w \in \R^7$ and $(v_{ij}) = (v_{ij} \ | \ 2 \leq i,j \leq 8)$.

By Lemma \ref{rklem} we know that $Y$ has rank $6$, that is, the column space of $Y$ is a six--dimensional subspace of $\R^7$.  The condition $[X, Y] = 0$ is equivalent to $w \in \R^7$ being perpendicular to each of the columns of $Y$.  Therefore $w$, and hence $X$, is uniquely determined up to scaling by $Y$.  Let $\psi\co \m \lra \p$ be the (unique) map assigning $X \in \p$ to $Y \in \m$ such that, by abusing the notation of (\ref{p}) and (\ref{m}), $X = \psi(Y)$ is given by
\beq
\label{psi}
w^t = \psi(v_1, \dots, v_7) = (v_7, -v_2, v_1, -v_4, v_3, v_6, -v_5)
\eeq
It is easy to check that $[\psi(Y), Y] = 0$, and $\psi$ is clearly a linear isomorphism.  Moreover, $\psi$ is $\G_2$-equivariant since $\Ad_g [X,Y] = [\Ad_g X, \Ad_g Y]$, for all $g \in \G_2$, and by uniqueness.  Therefore every horizontal zero--curvature plane $\sigma$ is determined by a pair $(X, Y) = (\psi(Y), Y)$, with $X \in \p$, $Y \in \m$.

Given $[X,Y] = 0$, the conditions in (\ref{EschLemConds}) imply that $\check \sigma_1$ has zero curvature if and only if $[(\Ad_{A^t} X)_\p, (\Ad_{A^t} Y)_\p] = 0$.  But $(\G,\K)$ is a rank one symmetric pair and thus $\check \sigma_1$ has zero curvature if and only if $(\Ad_{A^t} X)_\p, (\Ad_{A^t} Y)_\p$ are linearly dependent, that is, if and only if either
        $$(\Ad_{A^t} X)_\p = 0 \ \ {\rm or} \ \ (\Ad_{A^t} Y)_\p = 0 \ \ {\rm or} \ \ (\Ad_{A^t} X)_\p = s (\Ad_{A^t} Y)_\p$$
for some $s \in \R - \{0\}$.

Suppose that $(\Ad_{A^t} X)_\p = s (\Ad_{A^t} Y)_\p$, for some $s \in \R - \{0\}$.  Then $\Ad_{A^t} (X - sY) \perp \p$.  Since elements of $\p^\perp = \k = \mf{so}(7)$ have vanishing determinant, it follows from the discussion on ranks in the proof of Lemma \ref{rklem} that $\rank (X - s Y) = \rank \Ad_{A^t} (X - s Y) \leq 6$.  But $X$ and $Y$ describe a horizontal zero--curvature plane, so the vector $w \in \R^7$ determining $X$ is orthogonal to the columns of $sY$, which itself has rank $6$.  Hence $\rank (X - s Y) = 7$, which is a contradiction.   Therefore either $(\Ad_{A^t} X)_\p = 0$ or $(\Ad_{A^t} Y)_\p = 0$.

Before we continue, we recall some simple facts about $A \in \SO(8)$.  We may write $A$ as
$$
A = \left(\begin{array}{c|ccc}
a_{11}  & a_{12} & \cdots & a_{18}  \\ \hline
a_{2 1} & & & \\
\vdots & & \wt A & \\
a_{8 1} & & &
\end{array}\right)
$$
where $\wt A$ is a $7 \x 7$ matrix.  Since $A A^t = I$, one can easily derive the following identities:
\begin{align}
\label{A1}
a_{11}^2 + \sum |a_{1k}|^2 &= 1 \\
\label{A2}
a_{11} (a_{k1}) + \wt A (a_{1k})^t &= 0 \\
\label{A3}
(a_{k1})(a_{k1})^t + \wt A \wt A^t &= I_{7 \x 7}
\end{align}
where $(a_{1k})$ and $(a_{k 1})$ are row and column vectors respectively (with $2 \leq k \leq 8$).

Consider $X \in \p$ as in (\ref{XYvec}) and suppose that $(\Ad_{A^t} X)_\p = 0$.  A simple computation shows that this is equivalent to
$$
a_{11} w^t \wt A = (a_{k1})^t w (a_{1k}) = \<(a_{k1}),w \> (a_{1k}).
$$
If we multiply both sides by $(a_{1k})^t$, then the identity (\ref{A2}) above yields
\begin{align*}
-a_{11}^2 \<(a_{k1}),w \> &= -a_{11}^2 w^t (a_{k1})\\
&= \<(a_{k1}),w \> (a_{1k}) (a_{1k})^t \\
&= \<(a_{k1}),w \> \sum |a_{1k}|^2.
\end{align*}
By (\ref{A1}), this reduces to $(a_{k1})^t w = \<(a_{k1}),w \> = 0$.  Hence $a_{11} w^t \wt A = 0$, i.e. either $a_{11} = 0$ or $w^t \wt A = 0$.  However, if $w^t \wt A = 0$, then $\Ad_{A^t} X = 0$.  Since $X \neq 0$, this gives a contradiction.  Therefore the condition $(\Ad_{A^t} X)_\p = 0$ is satisfied if and only if $a_{11} = 0$ and $\<(a_{21}, \dots, a_{81})^t, w\> = 0$.  It is clear that the set of such $w$ is six--dimensional.  For reasons of dimension we will therefore always be able to find a pair $(X,Y)$ describing a zero--curvature plane, since we need only ensure that $\Ad_{A^t} X, \Ad_{A^t} Y \perp \sph^1$.

On the other hand, consider $Y \in \m$ as in (\ref{XYvec}) such that $(\Ad_{A^t} Y)_\p = 0$.  It is again simple to show that this is equivalent to
$$
(a_{k1})^t (v_{ij}) \wt A = 0.
$$
If we multiply both sides by $\wt A^t$, then (\ref{A3}) gives
\begin{align*}
0 &= (a_{k1})^t (v_{ij}) (I_{7 \x 7} - (a_{k1})(a_{k1})^t)\\
&= (a_{k1})^t (v_{ij}) - ((a_{k1})^t (v_{ij})(a_{k1}))(a_{k1})^t\\
&= (a_{k1})^t (v_{ij})
\end{align*}
where the last equality follows from $(a_{k1})^t (v_{ij})(a_{k1}) = 0$, since all entries on the diagonal of $\Ad_{A^t} Y \in \mf{so}(8)$ must be zero.

Let $X = \psi(Y) \in \p$ be given by $w \in \R^7$.  Now $(a_{k1})^t (v_{ij}) = 0$ implies that the vector $(a_{21}, \dots, a_{81})^t \in \R^7$ is perpendicular to each column of the rank $6$ matrix $(v_{ij})$.  Hence either $a_{k1} = 0$, for all $2 \leq k \leq 8$, or $w = (a_{21}, \dots, a_{81})^t$.

If $a_{k1} = 0$, for all $2 \leq k \leq 8$, then by (\ref{A1}), (\ref{A2}) and (\ref{A3}) $A \in \O(7) \subset \SO(8)$.  Conversely, given $A \in \O(7)$, it is clear that $\Ad_{A^t}$ preserves the orthogonal decomposition $\p \oplus \k$. Thus $\Ad_{A^t} X$, with $X \in \p$, will always be orthogonal to $\sph^1$ and, for dimension reasons, there will always be a $Y \in \m$ such that $\Ad_{A^t} Y \perp \sph^1$.  Therefore it is always possible to find a pair $(X,Y)$ spanning a horizontal zero--curvature plane.

Finally, in the case $w = (a_{21}, \dots, a_{81})^t$, there is, up to scaling, a unique pair $(X, Y)$ determined by the coordinates of $A$.  Clearly $X$ is determined by $w$ while, by (\ref{psi}), $Y$ is given by
$$
(v_1, \dots, v_7) = (a_{41}, -a_{31}, a_{61}, -a_{51}, -a_{81}, a_{71}, a_{21}).
$$
It is easy to see that the condition $\<X, \Ad_A E_{78}\>_0 = 0$ is given by the algebraic expression:
\beq
\label{AdX}
\sum_{k = 2}^8 (a_{18}a_{k7} - a_{17}a_{k8})a_{k1} = 0
\eeq
A straightforward computation shows that that condition $\<Y, \Ad_A E_{78}\>_0 = 0$ is given by
\beq
\label{AdY}
\begin{split}
%\nonumber
0 &= (m_{23} + m_{58} - m_{67}) a_{41}
- (m_{24} - m_{57} - m_{68}) a_{31} \\
%\nonumber
&{} \qquad + (m_{25} - m_{38} + m_{47}) a_{61}
- (m_{26} + m_{37} + m_{48}) a_{51}\\
%\label{AdY}
&{} \qquad - (m_{27} - m_{36} - m_{45}) a_{81}
+ (m_{28} + m_{35} - m_{46}) a_{71}\\
%\nonumber
&{} \qquad + (m_{34} + m_{56} - m_{78}) a_{21}
\end{split}
\eeq
where $m_{k\ell} = a_{k8}a_{\ell 7} - a_{k 7}a_{\ell 8}$.  It is clear that if both (\ref{AdX}) and (\ref{AdY}) hold, then there is a horizontal zero--curvature plane at $(A, I)$.

We have shown that a horizontal zero--curvature plane occurs at $(A,I)$ if and only if one of the following conditions is satisfied:
\begin{enumerate}
\item[(a)] $a_{11} = 0$
\item[(b)] $A \in \O(7) \subset \SO(8)$
\item[(c)] Equations (\ref{AdX}) and (\ref{AdY}) hold.
\end{enumerate}

The locus of such points is clearly lower dimensional.  Moreover, extending the $\UU$ action to an action by $\SO(3) \x \G_2$ increases the number of conditions which must be satisfied in order for a zero--curvature plane to be horizontal.  Theorem \ref{thmA}\ref{QPnew1} now follows immediately.

\section{Generalised Eschenburg spaces}
\label{GenEsch}

Consider the rank one symmetric pairs $(\G, \K) = (\U(n+1), \U(1)\U(n))$ and $(\K, \H) = (\U(1)\U(n), \U(1)\U(1)\U(n-1))$ where $n \geq 2$ and the inclusions $\K \hra \G$ and $\H \hra \K$ are given by
\begin{gather*}
(z,B) \lmt 
\begin{pmatrix} 
z & \\ & B 
\end{pmatrix}, 
\ \ z \in \U(1), \ \ B \in \U(n) \\
\tag*{and}
(z,w,C) \lmt 
\left(z, \begin{pmatrix} 
            w & \\ & C 
         \end{pmatrix} \right), 
\ \ z, w \in \U(1), \ \ C \in \U(n-1)
\end{gather*}
respectively.  Let $\g, \k$ and $\h$ be the Lie algebras of $\G, \K$ and $\H$ respectively.  Let $\< X, Y \>_0 = - \Re \tr (XY)$ be a bi-invariant metric on $\G$.  With respect to $\met_0$ we thus have the orthogonal decompositions $\g = \p \oplus \k$ and $\k = \m \oplus \h$, where:
        $$\p = \left\{\left(\begin{array}{c|ccc}
                                 0 & & -\bar x^t & \\
                            \hline & & & \\
                                 x & & 0 & \\
                                   & & & \end{array} \right) \ \Bigg| \ x = \bpm x_2 \\ \vdots \\ x_{n+1} \epm \in \C^n\right\} \ \ {\rm and}$$
        $$\m = \left\{\left(\begin{array}{c|c|ccc}
                                 0 & & & & \\
                        \hline     & 0 & & -\bar y^t & \\
                        \hline     & & & & \\
                                   & y & & 0 & \\
                                   & & & & \end{array} \right) \ \Bigg| \ y = \bpm y_3 \\ \vdots \\ y_{n+1} \epm \in \C^{n-1} \right\}$$
As in Examples \ref{Eschlem1} and \ref{Eschlem2} we define a left-invariant, right $\K$--invariant metric $\met_1$ on $\G$ by
        \beq
            \< X,Y\>_1 = \<X, \Phi_1(Y)\>_0
        \eeq
where $\Phi_1(Y) = Y_\p + \lambda_1 Y_\k$, $\lambda_1 \in (0,1)$, and a left-invariant, right $\H$--invariant metric $\met_2$ on $\G$ via
        \beq
        \label{met2}
            \<X,Y\>_2 = \<X,\Psi(Y)\>_1 = \<X,\Phi_2(Y)\>_0
        \eeq
where $\Phi_2(Y) = Y_\p + \lambda_1 Y_\m + \lambda_1 \lambda_2 Y_\h$, $\lambda_2 \in (0,1)$, and $\Psi(Y) = \Phi_1^{-1} \Phi_2(Y) = Y_\p + Y_\m + \lambda_2 Y_\h$.

Equip $\G \x \G$ with the left-invariant, right $(\K \x \H)$--invariant product metric $\met_1 \oplus \met_2$.

%\ms
Consider the subgroup $\L_{p,q} \subset \K \x \H$ defined by
        $$\L_{p,q} = \{(\diag(z^{p_1}, \dots, z^{p_{n+1}}), \diag(z^{q_1}, z^{q_2}, B)) \ | \ z \in \sph^1, B \in \U(n-1)\}$$
where $p_1, \dots, p_{n+1}, q_1, q_2 \in \Z$.  $\L_{p,q}$ acts on $\G$ via
        \begin{align*}
            \G &\lra \G\\
            A &\lmt \diag(z^{p_1}, \dots, z^{p_{n+1}}) A \diag(\bar z^{q_1}, \bar z^{q_2}, B^{-1})
        \end{align*}
where $z \in \U(1)$ and $B \in \U(n-1)$.  It is not difficult to show that this action is free if and only if
        \beq
        \label{freeGenEsch}
            (p_{\sigma(1)} - q_1, p_{\sigma(2)} - q_2) = 1 \ \ \textrm{for all} \ \ \sigma \in S_{n+1}.
        \eeq
We denote the resulting biquotients $\G \bq \L_{p,q}$ by $\E^{4n-1}_{p,q}$ and remark that $n = 2$ gives the usual Eschenburg spaces (see Eschenburg \cite[Section 41]{E1}).

Recall the canonical diffeomorphism
        $$\E^{4n-1}_{p,q} = \G \bq \L_{p,q} \cong \Delta \G \backslash \G \x \G / \L_{p,q}$$
given in (\ref{doubling}).  Now, since $\L_{p,q} \subset \K \x \H$, there is a metric on $\E^{4n-1}_{p,q}$ induced from the product metric on $\G \x \G$.

From (\ref{horizsp}) it is easy to show that the horizontal subspace at a point $(A, I) \in \G \x \G$ is given by
        \beq
        \label{horizGenEsch}
            \mcal{H}_A = \{(- \Phi_1^{-1} (\Ad_{A^*} W), \Phi_2^{-1} (W)) \ | \ W_{\mf{u}(n-1)} = 0, \<W, \Ad_A P - Q\>_0 = 0 \}
        \eeq
where $A^* = \bar A^t$, $P = \diag(i p_1, \dots, i p_{n+1})$, $Q = \diag(i q_1, i q_2, 0, \dots, 0)$, and $\h = \mf{u}(1) \oplus \u(1) \oplus \u(n-1)$ as before.

\begin{prop}
\label{0curv}
Suppose that
$$
\sigma = \Span\left\{\left(
- \Phi_1^{-1} \left(\Ad_{A^*} \wt X \right), \Phi_2^{-1} \left(\wt X \right)
\right),
\left(
- \Phi_1^{-1} \left(\Ad_{A^*} \wt Y \right), \Phi_2^{-1} \left(\wt Y \right)
\right) \right\}
$$
is a horizontal zero--curvature plane at $(A, I) \in \G \x \G$.  Then $X = \Phi_1^{-1}(\wt X)$ and $Y = \Phi_1^{-1}(\wt Y)$ can be written in one of the following forms:
\begin{enumerate}
\item \label{caseA}
$X \in \g$  and $Y = \diag(i, 0, \dots, 0)$
\item \label{caseB}
$X \in \p \oplus \h$ and $Y = \diag(0,i,0, \dots, 0)$
\item \label{caseC}
$X \in \p \oplus \h$ and $Y = \diag(i,i,0,\dots,0)$
\item \label{caseD}
            $$X = \left(\begin{array}{c|ccc}
                                 0 & & -\bar x^t & \\
                            \hline & & & \\
                                 x & & 0 & \\
                                   & & & \end{array} \right) \ \ \textrm{and} \ \ Y = \left(\begin{array}{c|c|ccc}
                                                                                                    i & & & & \\
                                                                                            \hline     & i\beta & & -\bar y^t & \\
                                                                                            \hline     & & & & \\
                                                                                                       & y & & 0 & \\
                                                                                                       & & & & \end{array} \right)
        $$
where $x_2 \neq 0$, $\beta = 1 - \sum_{j=3}^{n+1} |y_j|^2$, and $x_j = -i x_2 y_j$ for $j = 3, \dots, n+1$
\item \label{caseE}
            $$X = \left(\begin{array}{c|ccc}
                                 i\alpha & & -\bar x^t &\\
                            \hline & & &  \\
                                 x & & 0 & \\
                                   & & &  \end{array} \right) \ \ \textrm{and} \ \ Y = \left(\begin{array}{c|c|ccc}
                                                                                                    0 & & & & \\
                                                                                            \hline     & i\beta & & -\bar y^t & \\
                                                                                            \hline     & & & & \\
                                                                                                       & y & & 0 & \\
                                                                                                       & & & & \end{array} \right)
        $$
where $x = (0, x_3, \dots, x_{n+1})^t \neq 0 \in \C^n$ and $\sum_{j=3}^{n+1} x_j \bar y_j = 0$.
\end{enumerate}
\end{prop}

\begin{proof}
From the discussion in Section \ref{Biqs} we know that the projections $\vsig_i$, $i = 1, 2$, onto the first and second factor must be two-dimensional zero--curvature planes with respect to $\met_1$ and $\met_2$ respectively.

Consider
        $$\vsig_2 = \Span\left\{\Phi_2^{-1} \left(\wt X \right), \Phi_2^{-1} \left(\wt Y \right) \right\} = \Span \left\{\Psi^{-1} (X), \Psi^{-1} (Y) \right\}$$
where $\Psi = \Phi_1^{-1} \Phi_2$.  $\vsig_2$ has zero--curvature with respect to $\met_2$ if and only if the equalities in (\ref{met2conds}) hold, namely if and only if
        $$0 = [X, Y] = [X_\k, Y_\k] = [X_\p,Y_\p] = [X_\m,Y_\m] = [X_\h,Y_\h].$$
Since $(\G,\K)$ is a rank-one symmetric pair, $[X_\p, Y_\p] = 0$ if and only if $X_\p$ and $Y_\p$ are linearly dependent.  Without loss of generality we may assume that $Y_\p = 0$.  Similarly, $(\K,\H)$ being a rank-one symmetric pair implies that $[X_\m, Y_\m] = 0$ if and only if $X_\m$ and $Y_\m$ are linearly dependent.  Without loss of generality we may assume that either $X_\m = 0$ or $Y_\m = 0$.  Thus we have two possibilities:
        \begin{align}
        \label{case1}
                X &= X_\p + X_\m + X_\h \ \ {\rm and} \ \ Y = Y_\h \ \ \ {\rm or}\\
                X &= X_\p + X_\h \ \ {\rm and} \ \ Y = Y_\m + Y_\h.
        \label{case2}
        \end{align}
Since $\sigma$ is horizontal and $\Phi_1$ simply scales $\k = \m \oplus \h$ by $\lambda_1 \in (0,1)$, then we must have $X_{\u(n-1)} = Y_{\u(n-1)} = 0$, where $\h = \u(1) \oplus \u(1) \oplus \u(n-1)$.  Therefore in both cases above we have
        $$X_\h = \diag(ia, ib, 0, \dots, 0), \ \  Y_\h = \diag(ic, id, 0, \dots, 0), \ \ {\rm some} \ a,b,c,d \in \R.$$
Clearly $[X_\h, Y_\h] = 0$.  Then $[X_\k, Y_\k] = [X_\m, Y_\h] + [X_\h, Y_\m]$.  In the case of (\ref{case1}) our zero--curvature condition is thus $0 = [X_\m, Y_\h]$, while for case (\ref{case2}) we have $0 = [X_\h, Y_\m]$.

Consider general vectors $Z = \diag(i \alpha, i \beta, 0, \dots, 0) \in \h$ and $W \in \m$.  Then $0 = [Z, W]$ if and only if either $\beta = 0$ or $W = 0$.  Applying this to case (\ref{case1}) we find (after rescaling)
        \begin{align}
        \label{case3}
            X &= X_\p + X_\m + X_\h \ \ {\rm and} \ \ Y = \diag(i,0, \dots, 0) \ \ {\rm or}\\
        \label{case4}
            X &= X_\p + X_\h \ \ {\rm and} \ \ Y = Y_\h.
        \end{align}
On the other hand, case (\ref{case2}) yields the added possibility
        \beq
        \label{case5}
            X = \left(\begin{array}{c|ccc}
                                 i\alpha & & -\bar x^t & \\
                            \hline & & & \\
                                 x & & 0 & \\
                                   & & & \end{array} \right) \in \p \oplus \h \ \ {\rm and} \ \ Y = Y_\m + Y_\h.
        \eeq

As (\ref{case3}) is already of the form \ref{caseA} in the proposition, we concentrate on cases (\ref{case4}) and (\ref{case5}).

The only zero--curvature condition remaining to us is $[X,Y] = 0$.  Since $Y_\p = 0$, this is equivalent to $[X_\p, Y_\k] = 0$.  Consider the general vectors
        $$U = \left(\begin{array}{c|ccc}
                                 0 & & -\bar u^t & \\
                            \hline & & & \\
                                 u & & 0 & \\
                                   & & & \end{array} \right) \in \p \ \ {\rm and} \ \ V = \left(\begin{array}{c|c|ccc}
                                                                                                                    i\gamma & & & & \\
                                                                                                            \hline     & i\delta & & -\bar v^t & \\
                                                                                                            \hline     & & & & \\
                                                                                                                       & v & & 0 & \\
                                                                                                                       & & & & \end{array} \right) \in \k$$
where $u = (u_2, \dots, u_{n+1})^t \in \C^n$, $v = (v_3, \dots, v_{n+1})^t \in \C^{n-1}$, and $\gamma, \delta \in \R$.  Then
        \beq
        \label{nasty}
        \begin{split}
            &{} [U,V] = 0\\
             &\iff
                i(\gamma - \delta)u_2 + \sum_{j=3}^{n+1} u_j \bar v_j = 0 
                   \ \ {\rm and} \ \
                \gamma u_j = - i u_2 v_j, \  j = 3, \dots, n+1.
        \end{split}
        \eeq        
Suppose $u_2 = 0$.  Then (\ref{nasty}) becomes $\gamma u_j = 0$, $j = 3, \dots, n+1$, and $\sum_{j=3}^{n+1} u_j \bar v_j = 0$.  This is satisfied if and only if either
        \begin{align}
            \label{nasty1} &{}u_j = 0 \ \ \textrm{for all} \ \ j = 2, \dots, n+1,
 \ \ {\rm or}\\
            \label{nasty2} &{}u_2 = 0 \ \ {\rm and} \ \ \gamma = 0 \ \ {\rm and} \ \ \sum_{j=3}^{n+1} u_j \bar v_j = 0.
        \end{align}
On the other hand, if we assume $u_2 \neq 0$ then (\ref{nasty}) becomes
        \beq
        \label{nasty3}
        \begin{split}
                &{} u_2 \neq 0
                        \ \ {\rm and} \ \
                \gamma u_j = -i u_2 v_j, \  j = 3, \dots, n+1,
                        \ \ {\rm and} \\
                &{} \qquad\qquad \delta = \frac{\gamma}{|u_2|^2} \left( |u_2|^2 - \sum_{j=3}^{n+1} |u_j|^2 \right).
        \end{split}
        \eeq        
Now, if we apply conditions (\ref{nasty1}), (\ref{nasty2}) and (\ref{nasty3}) to case (\ref{case4}) we arrive at (after rescaling where appropriate)
        \begin{align}
            \label{case4a} &{} X = \diag(i a ,i b,0, \dots, 0) \ {\rm and} \ Y = \diag(i c,i d,0, \dots, 0) \ \ {\rm or} \\
            \label{case4b} &{} X \in \p \oplus \h \ \ {\rm and} \ \ Y = \diag(0,i,0, \dots, 0) \ \ {\rm or} \\
            \label{case4c} &{} X \in \p \oplus \h \ \ {\rm and} \ \ Y = \diag(i,i,0,\dots,0).
        \end{align}
Since $X$ and $Y$ must span a two-plane, it is clear that $\diag(i,0,\dots,0)$ must lie in the plane spanned by the $X$ and $Y$ given in (\ref{case4a}).  Therefore zero--curvature planes described by (\ref{case4}) fall into one of the classes given by \ref{caseA}, \ref{caseB} and \ref{caseC} of the proposition.

For case (\ref{case5}) conditions (\ref{nasty1}), (\ref{nasty2}) and (\ref{nasty3}) imply that $X$ and $Y$ must have one of the following forms (after rescaling):
        \begin{gather}
        \label{case5a}
            X = \diag(i, 0, \dots, 0) \ \ {\rm and} \ \ Y \in \k 
        \\ 
        \tag*{or}
            X = \left(\begin{array}{c|ccc}
                                 i\alpha & & -\bar x^t & \\
                            \hline  & & &  \\
                                  x & & 0 & \\
                                    & & & \end{array} \right) \ \ {\rm and} \ \             Y = \left(\begin{array}{c|c|ccc}
                                                                                                    0 & & & & \\
                                                                                            \hline     & i\beta & & -\bar y^t & \\
                                                                                            \hline     & & & & \\
                                                                                                       & y & & 0 & \\
                                                                                                       & & & & \end{array} \right)
        \end{gather}
where $x = (0, x_3, \dots, x_{n+1})^t \neq 0 \in \C^n$ and $\sum_{j=3}^{n+1} x_j \bar y_j = 0$, or finally
        \beq
        \label{case5c}
            X = \left(\begin{array}{c|ccc}
                                 i\alpha & & -\bar x^t & \\
                            \hline & & & \\
                                 x & & 0 & \\
                                   & & & \end{array} \right) \ \ {\rm and} \ \ Y = \left(\begin{array}{c|c|ccc}
                                                                                                    i & & & & \\
                                                                                            \hline     & i\beta & & -\bar y^t & \\
                                                                                            \hline     & & & & \\
                                                                                                       & y & & 0 & \\
                                                                                                       & & & & \end{array} \right)
        \eeq
where $x_2 \neq 0$, $\beta = 1 - \sum_{j=3}^{n+1} |y_j|^2$, and $x_j = -i x_2 y_j$ for $j = 3, \dots, n+1$.

Therefore, in order to complete the proof we may restrict our attention to horizontal zero--curvature planes for which $X$ and $Y$ are of the form (\ref{case5c}).

Without loss of generality we may assume that the vectors $\Psi^{-1}(X)$ and $\Psi^{-1}(Y)$ spanning $\vsig_2$ are orthogonal.  By (\ref{met2}) and since $Y \in \k$ this is equivalent to $\<X_\h, Y_\h\>_0 = 0$, where we recall that $\<V,W\>_0 = -\Re\tr(VW)$.  For (\ref{case5c}) we get orthogonality if and only if $\alpha = 0$.  Hence, as desired, we may rewrite (\ref{case5c}) as
        \Beq
         X = \left(\begin{array}{c|ccc}
                                 0 & & -\bar x^t & \\
                            \hline & & & \\
                                 x & & 0 & \\
                                   & & & \end{array} \right) \ \ {\rm and} \ \ Y = \left(\begin{array}{c|c|ccc}
                                                                                                    i & & & & \\
                                                                                            \hline     & i\beta & & -\bar y^t & \\
                                                                                            \hline     & & & & \\
                                                                                                       & y & & 0 & \\
                                                                                                       & & & & \end{array} \right)
        \Eeq
where $x_2 \neq 0$, $\beta = 1 - \sum_{j=3}^{n+1} |y_j|^2$, and $x_j = -i x_2 y_j$ for $j = 3, \dots, n+1$.
\end{proof}

In order to simplify the statements and computations to follow we fix
        $$A_t \ceq \left( \begin{array}{cc|ccc}
                        \cos t & - \sin t &&&\\
                        \sin t & \cos t &&& \\
                        \hline &&&& \\
                        &&& I_{n-1} & \\
                        &&&&
            \end{array} \right) \in \U(n+1), \ \ \ \ \cos^2 t \neq \Q$$
for the remainder of the proof, where $I_{n-1}$ denotes the $(n-1) \x (n-1)$ identity matrix.

\begin{lem}
\label{ABCnothoriz}
   If $p_1 \!\neq\! p_2$ and $p_1 + p_2 \neq q_1 + q_2$, then a vector $(- \Phi_1^{-1} (\Ad_{A^*} W), \Phi_2^{-1} (W))$ with
            $$W = \diag(i,0,\dots,0), \ \diag(0,i,0,\dots,0) \ \ {\rm or} \ \diag(i,i,0,\dots,0)$$
    cannot be horizontal at $(A_t, I) \in \G \x \G$.
\end{lem}

\begin{proof}
    Consider $V = \diag(i\theta, i\vphi, 0, \dots, 0)$.  From (\ref{horizGenEsch}) we see that a vector of the form $(- \Phi_1^{-1} (\Ad_{A_t^*} V), \Phi_2^{-1} (V))$ is horizontal if and only if $\<V, \Ad_{A_t} P - Q\>_0 = 0$.  Since $\< X, Y \>_0 = - \Re \tr (XY)$, this is equivalent to the condition
\begin{align*}
\theta q_1 + \vphi q_2 &= \theta \sum_{k=1}^{n+1} |a_{1k}|^2 p_k + \vphi \sum_{k=1}^{n+1} |a_{2k}|^2 p_k \\
&= (\theta \cos^2 t + \vphi \sin^2 t)p_1 + (\theta \sin^2 t + \vphi \cos^2 t)p_2.
\end{align*}
   If $(\theta, \vphi) = (1,0)$ this becomes $q_1 = p_1 \cos^2 t + p_2 \sin^2 t$ which is equivalent to $(p_1 - p_2)\cos^2 t = q_1 - p_2$, that is $\cos^2 t \in \Q$.  This is impossible by definition of $A_t$.  Similarly, when $(\theta, \vphi) = (0,1)$ we again have a contradiction.  Finally, when $(\theta, \vphi) = (1,1)$ we get $q_1 + q_2 = p_1 + p_2$, which contradicts our hypothesis.
\end{proof}

Recall that the projection $\vsig_1$ of a plane $\sigma$ as in Proposition \ref{0curv} onto the first factor is given by
        \begin{align*}
            \vsig_1 &= \Span\left\{\Phi_1^{-1} \left(\Ad_{A^*} \wt X \right), \Phi_1^{-1} \left(\Ad_{A^*} \wt Y \right) \right\}\\
                    &= \Span\left\{\Phi_1^{-1} \left(\Ad_{A^*} \left(\Phi_1 X \right) \right), \Phi_1^{-1} \left(\Ad_{A^*} \left(\Phi_1 Y \right)\right) \right\}.
        \end{align*}
(\ref{met1conds}) provides us with conditions for $\vsig_1$ to have zero--curvature with respect to $\met_1$ but, if we assume that $\vsig_2$ has zero--curvature, since $(\G,\K)$ is a rank-one symmetric pair the conditions reduce to
        \beq
        \label{sig1cond}
        \left[ \left( \Ad_{A^*} \left(\Phi_1 X \right) \right)_\p,
           \left( \Ad_{A^*} \left(\Phi_1 Y \right) \right)_\p \right] = 0.
        \eeq
That is, the $\p$ components of $\Ad_{A^*} \left(\Phi_1 X \right)_\p$ and $\Ad_{A^*} \left(\Phi_1 Y \right)_\p$ must be linearly dependent.  There are three possible cases:
        \begin{align*}
            &{} \left(\Ad_{A^*} \left(\Phi_1 X \right)\right)_\p = 0 \ \ {\rm or}\\
            &{} \left(\Ad_{A^*} \left(\Phi_1 Y \right)\right)_\p = 0 \ \ {\rm or}\\
             &{} \left(\Ad_{A^*} \left(\Phi_1 X \right)\right)_\p = s \left(\Ad_{A^*} \left(\Phi_1 Y \right)\right)_\p
        \end{align*}
for some $s \in \R-\{0\}$.

Recall that $\wt W = \Phi_1(W) = W_\p + \lambda_1 W_\k$, $\lambda_1 \in (0,1)$, and that
        $$\Ad_{A^*} \wt W = \sum_{k, \ell = 1}^{n+1} \bar a_{ki} a_{\ell j} \wt w_{k \ell}$$
where $\wt W = (\wt w_{ij})$ and $W \in \g$.  Then, since $\p$ is completely determined by the first row of vectors in $\g$, we may abuse notation to write
\beq
\label{AdAW}
(\Ad_{A^*} \wt W)_\p =
\left(\sum_{k, \ell = 1}^{n+1} \bar a_{k1} a_{\ell j} \wt w_{k \ell} \ \Big| \ j = 2, \dots, n+1 \right).
\eeq

\begin{lem}
\label{Dnothoriz}
Let $\sigma$ be a zero--curvature plane at $(A_t,I) \in \G \x \G$ as given by Proposition \ref{0curv}\ref{caseD}.  If $p_1 \neq p_2$ and $p_1 + p_2 \neq q_1 + q_2$ then $\sigma$ cannot be horizontal.
\end{lem}

\begin{proof}
By (\ref{AdAW}) we have
        \begin{align*}
            (\Ad_{A_t^*} (\Phi_1 X))_\p &= \left(-\cos^2 t \ \bar x_2 - \sin^2 t \ x_2 , - \cos t \ \bar x_3, \dots, - \cos t \ \bar x_{n+1} \right) \\
            (\Ad_{A_t^*} (\Phi_1 Y))_\p &= \left(i\lambda_1 (\beta - 1)\cos t \sin t, - \lambda_1 \sin t \ \bar y_3, \dots, - \lambda_1 \sin t \ \bar y_{n+1} \right).
        \end{align*}
If $(\Ad_{A_t^*} (\Phi_1 X))_\p = 0$ then $x_j = y_j = 0$ for all $j = 3, \dots, n+1$, since $x_j = -i x_2 y_j$, $j = 3, \dots, n+1$.  Thus $\beta = 1$ and $Y = \diag(i,i,0,\dots,0)$.  By Lemma \ref{ABCnothoriz} $Y$ cannot determine a horizontal vector.

If $(\Ad_{A_t^*} (\Phi_1 Y))_\p = 0$ then $Y = \diag(i,i,0,\dots,0)$ and, again, Lemma \ref{ABCnothoriz} shows that $\sigma$ cannot be horizontal.

Finally we examine the situation $\left(\Ad_{A_t^*} \left(\Phi_1 X \right)\right)_\p = s \left(\Ad_{A_t^*} \left(\Phi_1 Y \right)\right)_\p$ for some non-zero $s \in \R$.  Then $\cos t \ x_j = s \lambda_1 \sin t \ y_j$ implies that $-i \cos t \ x_2 y_j = s \lambda_1 \sin t \ y_j$, $j = 3, \dots, n+1$, since $x_j = -i x_2 y_j$, $j = 3, \dots, n+1$.  We have already shown that if $\sigma$ is to be horizontal the $y_j$ cannot all be zero.  Therefore $x_2 = i \left(\tfrac{s \lambda_1 \sin t}{\cos t}\right) \in i\R$.

Now $\left(-\Phi_1^{-1}\left(\Ad_{A_t^*} \left(\Phi_1 X \right) \right), \Phi_2^{-1}(\Phi_1 X) \right)$ is a horizontal vector if and only if the equation $\<X, \Ad_{A_t} P - Q\>_0 = 0$ is satisfied; that is, if and only if
        \begin{align*}
            0 &= \sum_{\ell = 1}^{n+1} \Im \left(a_{1\ell} \sum_{k=2}^{n+1} \bar a_{k\ell} x_k \right) p_\ell \\
            &= \Im \left(a_{11} \bar a_{21} x_2 \right) p_1 + \Im \left(a_{12} \bar a_{22} x_2 \right) p_2 \\
            &= \cos t \sin t \,(p_1 - p_2) \Im(x_2)
        \end{align*}
where again we recall that $\<V,W\>_0 = -\Re\tr(VW)$.  By hypothesis and definition of $A_t$ we must have $x_2 = i\Im(x_2) = 0$, which contradicts the assumption that $x_2 \neq 0$.
\end{proof}

\begin{lem}
\label{Enothoriz}
Let $\sigma$ be a plane at $(A_t, I) \in \G \x \G$ as determined by Proposition \ref{0curv}\ref{caseE}.  Then $\sigma$ does not have zero curvature at $(A_t,I)$.
\end{lem}

\begin{proof}
Following (\ref{AdAW}) we write
        \begin{align*}
            (\Ad_{A_t^*} (\Phi_1 X))_\p &= \left(-i \lambda_1 \alpha \cos t \sin t, - \cos t \ \bar x_3, \dots, - \cos t \ \bar x_{n+1} \right) \\
            (\Ad_{A_t^*} (\Phi_1 Y))_\p &= \left(i \lambda_1 \beta \cos t \sin t, - \lambda_1 \sin t \ \bar y_3, \dots, - \lambda_1 \sin t \ \bar y_{n+1} \right).
        \end{align*}
If $(\Ad_{A_t^*} (\Phi_1 X))_\p = 0$ then $X = 0$, which is a contradiction since $\sigma$ is two-dimensional.  Similarly, $(\Ad_{A_t^*} (\Phi_1 Y))_\p = 0$ gives a contradiction.  On the other hand, if there is some non-zero $s \in \R$ such that $\left(\Ad_{A_t^*} \left(\Phi_1 X \right)\right)_\p = s \left(\Ad_{A_t^*} \left(\Phi_1 Y \right)\right)_\p$ then we find $\cos t \ x_j = s \lambda_1 \sin t \ y_j$, $j = 3, \dots, n+1$.  However, since $\sum_{j=3}^{n+1} x_j \bar y_j = 0$, this implies that $x_j = y_j = 0$ for all $j = 3, \dots, n+1$.  But, by hypothesis, $x_j$ cannot all be zero for planes of this type, and so we have a contradiction.  Therefore, by the discussion following (\ref{sig1cond}), $\sigma$ cannot admit zero--curvature at $(A_t, I)$.
\end{proof}

\pf{\ref{thmB}}
If $p_1 = \dots = p_{n+1}$ then we may assume without loss of generality that $p_i = 0$, for all $i = 1, \dots, n+1$.  Since the action is free we must have $(q_1, q_2) \neq (0,0)$.  Hence, by Tapp \cite{Ta1}, we are done.

As discussed in Section \ref{Biqs}, a permutation of the integers $p_1, \dots, p_{n+1}$ induces a diffeomorphism $\E^{4n-1}_{p,q} \lra \E^{4n-1}_{p,q}$.  Therefore, if $p_i \neq p_j$ for some $i \neq j$, we may assume that $p_1 \neq p_2$.  If $p_1 + p_2 \neq q_1 + q_2$ then we are done by Proposition \ref{0curv} and Lemmas \ref{ABCnothoriz}, \ref{Dnothoriz} and \ref{Enothoriz}.

Finally, if $p_1 \neq p_2$ but $p_1 + p_2 = q_1 + q_2$, then the freeness condition (\ref{freeGenEsch}) implies that either $(p_1, p_2, \dots, p_{n+1}) = (1,\dots, 1, -1,\dots, -1)$ and $(q_1, q_2) = (0,0)$ (up to sign and permutations of the $p_i$), where $\#\{i \ | \ p_i = -1\} \in \left\{\left\lfloor \tfrac{n+2}{2}\right\rfloor, n\right\}$, or there is some $p_k \not\in \{p_1, p_2\}$, $k \in \{3, \dots, n+1\}$.

In the latter case we may permute and relabel the $p_i$ such that $p_1 + p_2 \neq q_1 + q_2$, in which case we are done as above.  In the former case we may assume that $p_1 = 1$, $p_2 = -1$ and $p_3 = -1$ (after reordering if necessary).  Then, for
$$
A_0 \ceq \left( \begin{array}{ccc|ccc}
                        \frac{1}{2} & \frac{1}{\sqrt{2}} & \frac{1}{2} &&&\\
                        -\frac{1}{2} & \frac{1}{\sqrt{2}} & -\frac{1}{2} &&& \\
                        -\frac{1}{\sqrt{2}} & 0 & \frac{1}{\sqrt{2}} &&& \\
                        \hline &&&&& \\
                        &&&& I_{n-2} & \\
                        &&&&&
            \end{array} \right) \in \U(n+1),$$
we may repeat the approaches used in the proofs of Lemmas \ref{ABCnothoriz}, \ref{Dnothoriz} and \ref{Enothoriz} to rule out the existence of any horizontal zero--curvature planes at $(A_0, I) \in \G \x \G$.  Note that, unlike in the proof of Lemma \ref{Enothoriz}, it is necessary to show that a plane as determined by Proposition \ref{0curv}\ref{caseE} cannot have zero curvature and be horizontal at the same time.
\qed

\section{Topology of $\M^{13}$ and $\N^{11}$}
\label{TopMN}

We turn now to the topological assertions of Theorem \ref{thmA} regarding the biquotients $\M^{13} = \SO(8) \bq (\sph^1 \x \G_2)$ and $\N^{11} = \SO(8) \bq (\SO(3) \x \G_2)$, namely that they have the same cohomology rings but are not homeomorphic to $\CP^3 \x \sph^7$ and $\sph^4 \x \sph^7$ respectively.

\begin{thm}
\label{cohom}
The biquotients $\M^{13}$ and $\N^{11}$ have the same cohomology rings as $\CP^3 \x \sph^7$ and $\sph^4 \x \sph^7$ respectively.  In particular $\M^{13}$ and $\N^{11}$ are not manifolds known to admit positive curvature.
\end{thm}

\begin{proof}
In Section \ref{SOFree} we established that $\M^{13}$ and $\N^{11}$ are the total spaces of $\sph^7$--bundles over $\CP^3$ and $\sph^4$ respectively.  Given an arbitrary fibration $\sph^7 \lra \E \lra \B$, where $\B$ is a simply connected, compact manifold with $\dim \B \leq 7$, the Euler class $e$ is trivial since $\H^8(\B; \Z) = 0$ and so the Gysin sequence splits into short exact sequences, from which it immediately follows that $\H^j (\E; \Z) = \H^j (\B \x \sph^7; \Z)$, for all $j$.  Now $\M^{13}$ and $\N^{11}$ are quotients of $\sph^7 \x \sph^7$ by $\sph^1$ and $\sph^3$ respectively.  A quick glance at the Serre spectral sequences of these fibrations shows that the ring structures of $\H^* (\M^{13}; \Z)$ and $\H^* (\N^{11}; \Z)$ agree with those of $\H^* (\CP^3 \x \sph^7; \Z)$ and $\H^* (\sph^4 \x \sph^7; \Z)$ as desired.
\end{proof}

In order to distinguish $\M^{13}$ and $\N^{11}$ from the products $\CP^3 \x \sph^7$ and $\sph^4 \x \sph^7$, respectively, we want to analyse their Pontrjagin classes.  There is a general procedure (developed in Borel and Hirzebruch \cite{BH}, Eschenburg \cite{E2} and Singhof \cite{Si}) for computing the $\RR$--Pontrjagin class of a biquotient $\G \bq \UU$, where $\RR$ is a coefficient ring such that $\H^*(\G; \RR)$ and $\H^*(\UU; \RR)$ have no torsion, and the action of $\UU$ on $\G$ is (effectively) free.  Let $\iota \co \UU \lra \GG$ denote the embedding and assume that we have such an $\RR$.

We adopt the following notation: For a compact, connected Lie group $\L$, let $\T_\L$ denote the maximal torus and $W_\L$ the Weyl group.  Let $\E_\L$ be a contractible space on which $\L$ acts freely.  The classifying space of $\L$ is the quotient $\B_\L \ceq \E_\L/\L$.  A product of Lie groups $\L_1$ and $\L_2$ is written $\L_1 \L_2$.

Consider the following commutative diagram of fibrations
        $$
            \xymatrix{ \G \x \E_{\GG} \ar[r] \ar[d] & \G \x \E_{\GG} \ar[d] & \\
                        \G \x_{\UU} \E_{\GG} \ar[r]^{\varphi_{\G}} \ar[d]_{\varphi_\UU} & \G \x_{\GG} \E_{\GG} \ar[d]^{\B_\Delta} & \hspace{-.9cm} = \B_{\Delta \G} \\
                        \B_\UU \ar[r]_{\B_\iota} & \B_{\GG} &
        }$$
where $\varphi_{\G}$ and $\varphi_\UU$ are the respective classifying maps, and $\Delta \co \G \lra \GG$ denotes the diagonal embedding.  Now, since projection onto the first factor in each case is a homotopy equivalence, we have $\G \simeq \G \x \E_{\GG}$ and $ \G \bq \UU \simeq \G \x_{\UU} \E_{\GG}$.  Thus, up to homotopy, we can consider the diagram above as
        \beq
        \label{diagram}
           \xymatrix{ \G \ar[r] \ar[d] & \G \ar[d] \\
                      \G \bq \UU \ar[r]^{\varphi_{\G}} \ar[d]_{\varphi_\UU} &
                         \B_{\G} \ar[d]^{\B_\Delta} \\
                      \B_\UU \ar[r]_{\B_\iota} & \B_{\GG}
        }\eeq
We would like to use the Serre spectral sequences of the fibrations on the left and right, as well as the commutativity of the diagram, to compute the Pontrjagin class of (the tangent bundle of) $\G \bq \UU$.  We will need to know $\H^*(\G; \RR)$, $\H^*(\B_\UU; \RR)$, $\H^*(\B_{\GG}; \RR)$ and the homomorphisms $(\B_\iota)^* \co \H^*(\B_{\GG}; \RR) \lra \H^*(\B_\UU; \RR)$ and $(\varphi_\UU)^* \co \H^*(\B_\UU; \RR) \lra \H^*(\G \bq \UU; \RR)$.

Let $\L$ be a compact, connected Lie group and $\RR$ a ring such that $\H^*(\L; \RR)$ has no torsion.  Then
        \beq
        \label{cohomL}
            \H^*(\L; \RR) = \Lambda (y_1, \dots, y_r)
        \eeq
where $r = \rank(\L)$.  Hence, from the Serre spectral sequence of the universal bundle $\L \lra \E_\L \lra \B_\L$,
        \beq
            \label{cohomBL}
            \H^*(\B_\L; \RR) = \RR [\bar y_1, \dots, \bar y_r]
        \eeq
where $\bar y_j$ denotes the transgression of $y_j$.

Let $(t_1, \dots, t_r)$ be coordinates of the maximal torus $\T_\L$.  By an abuse of notation we identify $t_j$ with the element $t_j \in \H^1(\T_\L; \RR)$.  The corresponding transgression arising from the Serre spectral sequence for $\T_\L \lra \E_{\T_\L} \lra \B_{\T_\L}$ is $\bar t_j \in \H^2 (\B_{\T_{\G}}; \RR)$.  Since $\L$ does not have any torsion in its $\RR$-cohomology we have
        \beq
           \label{cohomGENS}
           \H^*(\B_\L; \RR) =
              \H^*(\B_{\T_\L}; \RR)^{W_\L} = \RR [\bar t_1, \dots, \bar t_r]^{W_\L}
        \eeq
Thus we have explicit well-defined generators of $\H^*(\B_\L; \RR)$ which we identify with $\bar y_j$, $j = 1, \dots, r$.

Suppose now that $h \co \L_1 \lra \L_2$ is a homomorphism of Lie groups.  Then the commutative diagram
        $$\xymatrix{\L_1 \ar[r]^h & \L_2 \\ \T_{\L_1} \ar@{^{(}->}[u] \ar[r]_h & \T_{\L_2} \ar@{^{(}->}[u] }$$
induces a commutative diagram of classifying spaces
        \beq \label{classdiagram}
            \xymatrix{\B_{\L_1} \ar[r]^{\B_h} & \B_{\L_2} \\ \B_{\T_{\L_1}} \ar[u] \ar[r]_{\B_h} & \B_{\T_{\L_2}} \ar[u] }
        \eeq
which in turn induces the commutative diagram
        \beq \label{cohomdiagram}
            \xymatrix{\H^*(\B_{\L_1}) \ar[d] &  \H^*(\B_{\L_2}) \ar[l]_{(\B_h)^*} \ar[d] \\ \H^*(\B_{\T_{\L_1}})  & \H^*(\B_{\T_{\L_2}}) \ar[l]^{(\B_h)^*} }
        \eeq
In particular, by (\ref{cohomGENS}), we see that $(\B_h)^*$ can be determined by simply understanding $h(\T_{\L_1}) \subset \T_{\L_2}$.

Consider again diagram (\ref{diagram}).  Recall that, since there is no torsion,
        \Beq
            \H^*(\B_{\GG}; \RR) = \H^*(\B_{\G}; \RR) \ox \H^*(\B_{\G}; \RR).
        \Eeq
Then $\H^*(\B_{\GG}; \RR)$ is generated by classes of the form $\bar y_j \ox 1$ and $1 \ox \bar y_j$, $j = 1, \dots, r= \rank(\G)$.  Consider the diagonal embedding $\Delta \co \G \hra \GG$.  In coordinates $\Delta|_{\T_{\G}}$ is given by $t_j \lmt (t_j, t_j)$, $j=1, \dots, r$.  We have commutative diagrams as in (\ref{classdiagram}) and (\ref{cohomdiagram}).  Now
        \begin{align*}
            \Delta^* \co \H^1(\T_{\G}; \RR) \ox \H^1(\T_{\G}; \RR)  &\lra \H^1(\T_{\G}; \RR)\\
                        t_j \ox 1 &\lmt t_j\\
                        1 \ox t_j &\lmt t_j
        \end{align*}
which in turn implies
        \begin{align*}
            (\B_\Delta)^* \co \H^2(\B_{\T_{\G}}; \RR) \ox \H^2(\B_{\T_{\G}}; \RR) &\lra \H^2(\B_{\T_{\G}}; \RR)\\
                            \bar t_j \ox 1 &\lmt \bar t_j\\
                            1 \ox \bar t_j &\lmt \bar t_j.
        \end{align*}
Therefore, by (\ref{cohomGENS}),
        \begin{align*}
            (\B_\Delta)^* \co \H^*(\B_{\GG}; \RR) &\lra \H^*(\B_{\G}; \RR)\\
                        \bar y_j \ox 1 &\lmt \bar y_j\\
                        1 \ox \bar y_j &\lmt \bar y_j.
        \end{align*}
Since the diagram (\ref{diagram}) is commutative we see that
        \begin{align}
            \nonumber
            \vphi_{\G}^*(\bar y_j) &= \vphi_{\G}^* ((\B_\Delta)^*(\bar y_j \ox 1)) = \vphi_\UU^*((\B_{\iota})^*(\bar y_j \ox 1)) \ \ \ {\rm and}\\
            \label{equality}
            \vphi_{\G}^*(\bar y_j) &= \vphi_{\G}^* ((\B_\Delta)^*(1 \ox \bar y_j)) = \vphi_\UU^*((\B_\iota)^*(1 \ox \bar y_j)).
        \end{align}
We remark that, by naturality of spectral sequences, one can deduce from (\ref{diagram}) that
$$
\H^* (\G \bq \UU; \RR) = \H^* \left(\frac{\H^*(\G; \RR) \ox \H^*(\B_\UU; \RR)}{\left\<dy_j = \vphi_\UU^*((\B_{\iota})^*(\bar y_j \ox 1)) - \vphi_\UU^*((\B_\iota)^*(1 \ox \bar y_j)) \right\>}\right).
$$

Let us now focus on computing the Pontrjagin class of $\G \bq \UU$.  Let $\tau$ be the tangent bundle of $\G \bq \UU$.  In analogue with Singhof \cite{Si}, we introduce the following vector bundles over $\G \bq \UU$:  First, let $\alpha_{\G} \ceq (\GG / \UU) \x_{\G} \g$, where $\G$ acts on $(\GG / \UU) \x \g$ via
$$
g \star ([g_1, g_2], X) = ([g g_1 , g g_2], \Ad_g X)
$$
with $g, g_1, g_2 \in \G$, $[g_1, g_2] = (g_1 , g_2) \cdot \UU \in \GG/\UU$, and $X \in \g$.  Second, let $\alpha_\UU \ceq \G \x_\UU \u$, where $\u \subset \g \oplus \g$ is the Lie algebra of $\UU$ and $\UU$ acts on $\G \x \u$ via
$$
(u_1, u_2) \star (g, (Y_1, Y_2)) = (u_1 g u_2^{-1}, (\Ad_{u_1} Y_1, \Ad_{u_2} Y_2))
$$
with $g \in \G$, $(u_1, u_2) \in \UU$ and $(Y_1, Y_2) \in \u$.  Then, as in Singhof \cite[Proposition 3.2]{Si}:

\begin{lem}
\label{directsum}
$\tau \oplus \alpha_\UU = \alpha_{\G}$.
\end{lem}

\begin{proof}
Recall that, since $\UU$ acts on $\G$ via $(u_1, u_2) \star g = u_1 g u_2 ^{-1}$, $(u_1, u_2) \in \UU$, $g \in \G$, the vertical subspace at $g \in \G$ (after left-translating back to $e \in \G$) is given by
$$
\mcal{V}_g = \{\Ad_{g^{-1}} Y_1 - Y_2 \ | \ (Y_1, Y_2) \in \u\} \subset \g.
$$
Moreover, since the action of $\UU$ on $\G$ is free, we have
\beq
\label{freeInt}
\{(Y_1, Y_2) \in \u \ | \ \Ad_{g^{-1}} Y_1 = Y_2 \ \textrm{for some} \ g \in \G\} = \{0\}.
\eeq
Given $g \in \G$, consider the maps
\begin{gather*}
\begin{aligned}
\psi_g \co \g &\lra \alpha_{\G} \\
X &\lmt [[g^{-1},e], X]
\end{aligned}\\
\tag*{and}
\begin{aligned}
f \co \alpha_\UU &\lra \alpha_{\G} \\
[g, (Y_1, Y_2)] &\lmt [[g^{-1},e], \Ad_{g^{-1}} Y_1 - Y_2].
\end{aligned}
\end{gather*}
It is easy to check that $f$ is well-defined and that $f_{[g]}((\alpha_\UU)_{[g]}) = \psi_g (\mcal{V}_g)$, where $(\alpha_\UU)_{[g]}$ is the fibre of the bundle $\alpha_\UU$ over $[g] = g \cdot \UU \in \G \bq \UU$.  Furthermore, (\ref{freeInt}) allows us to deduce that $f$ is injective.

If we now equip $\G$ with a bi-invariant metric, then (after left-translation to the identity in $\G$) the horizontal subspace at $g$ is given by
$$
\mcal{H}_{g} = \left\{ X \ | \ X \in \g, \  X \perp \mcal{V}_{g} \right\} \subset \g.
$$
Note that $\UU$ acts on the horizontal distribution by taking $X \in \mcal{H}_{g}$ to $\Ad_{u_2} X \in \mcal{H}_{u_1 g u_2^{-1}}$, where $(u_1, u_2) \in \UU$.  Let $\check{\mcal{H}} = \mcal{H}/\UU$ be the orbit space under the $\UU$ action.  In particular, $\tau \cong \check{\mcal{H}}$ and we can represent the tangent space at $[g] \in \G \bq \UU$ by
$$
\check{\mcal{H}}_{[g]} = \left\{ [g, X] \ | \ X \in \g, \  X \perp \mcal{V}_g \right\}
$$
where the $g$ in $[g, X]$ serves only to keep track of the base point.

Recall that each point in $\alpha_{\G}$ may be represented as $[[g_1, g_2], \Ad_{g_2} X] = [[g_2^{-1} g_1, e], X]$, for some $g_1, g_2 \in \G$ and $X \in \g$.  Consider the map
\begin{align*}
q \co \alpha_{\G} &\lra \check{\mcal{H}}_{[g]} = \tau \\
[[g_1,g_2], \Ad_{g_2} X] &\lmt [g_1^{-1} g_2, X^\perp ]
\end{align*}
where $X^\perp$ is the component of $X \in \g$ orthogonal to $\mcal{V}_{g_1^{-1} g_2}$.  Clearly $q$ is surjective.  The fact that $q$ is well-defined follows from noticing that $(\Ad_{u_2} X)^\perp = \Ad_{u_2} (X^\perp) \in \mcal{H}_{u_1 g u_2^{-1}}$, for all $X \in \g$, $g \in \G$ and $(u_1, u_2) \in \UU$.

Now $q([[g_1,g_2], \Ad_{g_2} X]) = 0$ for some $[[g_1,g_2], \Ad_{g_2} X] \in \alpha_{\G}$ if and only if $X^\perp = 0 \in \mcal{H}_{g_1^{-1} g_2}$, i.e. $X \in \mcal{V}_{g_1^{-1} g_2}$.  Thus $[[g_1,g_2], \Ad_{g_2} X] = [[g_2^{-1} g_1, e], X] \in f_{[g_1^{-1} g_2]}\left((\alpha_\UU)_{[g_1^{-1} g_2]}\right)$ and we have $\ker q = f(\alpha_\UU)$.  Therefore $\alpha_{\G} / f(\alpha_\UU) = \tau$.  In other words, $\alpha_{\G} = \tau \oplus \alpha_\UU$ as desired.
\end{proof}

Recall from Borel--Hirzebruch \cite{BH} that, if $V$ is a representation of a Lie group $\L$, then the homogeneous vector bundle $\alpha_\L = \P \x_\L V$ associated to the $\L$--principal bundle $\P \lra \B\ceq \P/\L$ is the pullback under the classifying map $\varphi_\L \co \B \lra \B_\L$ of the vector bundle $\E_\L \x_\L V$ associated to the universal $\L$--principal bundle $\E_\L \lra \B_\L$, i.e. $\alpha_\L = \varphi_\L^*(\E_\L \x_\L V)$.  In particular, from Borel--Hirzebruch \cite{BH} and Singhof \cite{Si} we know that the Pontrjagin class of the bundle $\alpha_\L \lra \B$ is given by
        $$\pp(\alpha_\L) = 1 + \pp_1 (\alpha_\L) + \pp_2 (\alpha_\L) + \dots = \varphi_\L^* (a), \ \ \ a \ceq \prod_{\alpha_i \in \Delta_\L^+} (1 + \bar \alpha_i^2) $$
where $\Delta_\L^+$ is the set of positive weights of the representation of $\L$ on $V$.  We have identified $\alpha_i \in \H^1(\T_\L; \RR) \cong \H^2 (\B_{\T_\L}; \RR)$, and so it follows that $a \in \H^*(\B_{\T_\L}; \RR)^{W_\L} \cong \H^*(\B_\L; \RR)$.

The vector bundles $\alpha_\UU$ and $\alpha_{\G}$ are associated to principal $\UU$ and $\G$ bundles respectively and, in this case, the weights of the defining representations are the roots of the corresponding Lie groups.

Since $\pp(V \oplus W) = \pp(V) \smallsmile \pp(W)$, for vector bundles $V$ and $W$ over some manifold $\M$, we have
\vspace{.5mm}
        $$\pp(\tau)\pp(\alpha_\UU) = \pp(\alpha_{\G}).$$
By our discussion above and since inverses are well-defined in the polynomial algebra $\H^*(\B_\UU; \RR)$ it follows that
        $$\pp(\tau) = \varphi_{\G}^*(a) \varphi_\UU^*(b^{-1})$$
where $a \ceq \prod_{\alpha_i \in \Delta_{\G}^+} (1+\bar \alpha_i^2)$ and $b \ceq \prod_{\beta_j \in \Delta_\UU^+} (1+\bar \beta_j^2)$.  In particular,
        \begin{align}
           \nonumber
            \pp_1(\tau) &= \pp_1(\alpha_{\G}) - \pp_1(\alpha_\UU)\\
           \label{firstPont}
                      &= \varphi_{\G}^*\left( \sum_{\alpha_i \in \Delta_{\G}^+} \bar \alpha_i^2 \right) -  \varphi_\UU^* \left( \sum_{\beta_j \in \Delta_\UU^+} \bar \beta_j^2 \right).
        \end{align}

%\ms

In our situation, the Lie groups $\G = \SO(8)$ and $\UU_k \ceq \SO(k) \x \G_2$, with $k = 2,3$, have no torsion in their cohomology for coefficients in $\RR = \Z_p$, $p$ an odd prime (see Mimura and Toda \cite[Corollary 3.15, Theorem 5.12]{MT}).  We use the process outlined above to compute the $\Z_p$--Pontrjagin classes of $\M^{13}$ and $\N^{11}$.  Before we continue we establish an easy lemma which will prove useful in the topological computations to follow.

\begin{lem}
\label{nicetrick}
Consider a triple $(r_1, r_2, r_3)$ such that $\sum r_i = 0$.  Let $\sigma_i (r)$ and $\sigma_i (r^2)$ denote the $i^{\it th}$ elementary symmetric polynomials in $r_1, r_2, r_3$ and $r_1^2, r_2^2, r_3^2$ respectively.  Then $\sigma_1 (r^2) = -2 \sigma_2 (r)$ and $\sigma_2 (r^2) = \sigma_2 (r)^2$.
\end{lem}

\proof
Since $\sigma_1 (r) = \sum r_i = 0$ we have
        \begin{align*}
            0 &= \sigma_1 (r)^2\\
            &= (r_1^2 + r_2^2 + r_3^2) + 2(r_1 r_2 + r_1 r_3 + r_2 r_3)\\
            &= \sigma_1(r^2) + 2 \sigma_2(r)
        \end{align*}
as desired.  On the other hand:
        \begin{align*}
            \sigma_2(r)^2 - \sigma_2(r^2) &= (r_1 r_2 + r_1 r_3 + r_2 r_3)^2 - (r_1^2 r_2^2 + r_1^2 r_3^2 + r_2^2 r_3^2)\\
                                        &= 2(r_1^2 r_2 r_3 + r_1 r_2^2 r_3 + r_1 r_2 r_3^2)\\
                                        &= 2 r_1 r_2 r_3 (r_1 + r_2 + r_3)\\
                                        &= 0
            \rlap{\hspace{7.78cm} \qedsymbol}
        \end{align*}

In \cite[pages $vii$ and 139]{E1}, Eschenburg provides a beautiful diagram which explicitly describes the embedding of the root system $\G_2$ into $\B_3$.  Recall that $\B_3$ is the root system corresponding to the Lie algebra $\mf{so}(7)$ and is given by
        $$\B_3 = \{\pm t_i \ | \ 1 \leq i \leq 3\} \cup \{\pm (t_i \pm t_j) \ | \ 1 \leq i < j \leq 3\}.$$
The root system $\G_2$ lies on a hypersurface in $\Span \{\B_3\}$ and is given by
        $$\G_2 = \{\pm s_i \ | \ 1 \leq i \leq 3\} \cup \{\pm (s_i - s_j) \ | \ 1 \leq i < j \leq 3\}$$
where $s_i = \frac{1}{3}(2 t_i - t_j - t_k)$, $\{i,j,k\} = \{1,2,3\}$.  Notice that $\sum s_i = 0$ and that $s_i - s_j = t_i - t_j \in \B_3$.  Furthermore, $s_i$ is the projection of $t_i \in \B_3$ and $-(t_j + t_k) \in \B_3$ onto the hypersurface containing $\G_2$.

Since the Lie group $\G_2$ is simply connected and has no centre, we see that the inclusions
        $$\exp^{-1}(I) = \textrm{integral lattice of } \G_2 \subset \textrm{root lattice of } \G_2 \subset \textrm{weight lattice of } \G_2$$
are in fact equalities.  Therefore, by our above discussion of the roots of $\G_2$, the integral and weight lattices of $\G_2$ are spanned by $\{s_i \ | \ 1 \leq i \leq 3\}$, $\sum s_i = 0$.  Thus by an abuse of notation we may assume that $\{s_i \ | \ 1 \leq i \leq 3\}$, $\sum s_i = 0$, spans $\H^1 (\T_{\G_2}; \Z) = \Hom(\Gamma, \Z)$, where $\T_{\G_2}$ is a maximal torus of $\G_2$ and $\Gamma$ is the integral lattice of $\G_2$.

\begin{prop}
\label{Pontclass13}
    Let $p$ be an odd prime.  The $\Z_p$ first Pontrjagin classes of $\M^{13}$ and $\N^{11}$ are
        $$
            \pp_1 (\M^{13}) = 2 \alpha^2 \ \ \ \textrm{and} \ \ \ \pp_1 (\N^{11}) = \beta        
        $$
    respectively, where $\alpha$ is a generator of $\H^2(\M^{13}; \Z_p) = \Z_p$ and $\beta$ is a generator of $\H^4 (\N^{11}; \Z_p) = \Z_p$.
\end{prop}

%\pfP{\ref{Pontclass13}}
\begin{proof}
Let $\G = \SO(8)$ and let $\UU_k \ceq \SO(k) \x \G_2$, $k = 2, 3$, act freely on $\G$ as described in Section \ref{SOFree} with quotients $\M^{13}$ and $\N^{11}$ respectively. Let $\iota_k \co \UU_k \hra \GG$, $k = 2,3$, denote the respective inclusions.

For $p$ an odd prime, the $\Z_p$-cohomology of $\G$ and $\UU_k$ is
        \begin{align*}
            \H^*(\G; \Z_p) &= \Lambda (y_1, y_2, y_3, y_4), \ \ \ y_1 \in \H^3, y_2, y_4 \in \H^7, y_3 \in \H^{11} \\
            \H^*(\UU_k; \Z_p) &= \Lambda (w_k, x_1, x_2), \ \ \ w_k \in \H^{2k-3}(\SO(k)), \ x_i \in \H^{8i-5}(\G_2).
        \end{align*}

Let $\T_{\G}$ and $\T_\UU$ be the maximal tori of $\G$ and $\UU_k$ respectively, with coordinates being given by $(t_1, t_2, t_3, t_4)$ and $(u, s_1, s_2, s_3)$, $\sum s_i = 0$, respectively.  By an abuse of notation (and our earlier discussion of the roots of $\G_2$) we will identify $t_i$, $u$ and $s_j$ with the elements $t_i \in \H^1(\T_{\G})$ and $u, s_j \in \H^1 (\T_\UU)$.  The corresponding transgressions are $\bar t_i \in \H^2 (\B_{\T_{\G}})$ and $\bar s_j \in \H^2 (\B_{\T_\UU})$.

$W_{\G}$ acts on $\H^*(\B_{\T_{\G}})$ via permutations in $\bar t_i$ and an even number of sign changes.  Therefore a basis for $\H^*(\B_{\T_{\G}})^{W_{\G}}$ is given by elementary symmetric polynomials $\sigma_i (\bar t^2) \ceq \sigma_i (\bar t_1^2, \dots, \bar t_4^2)$, $i=1,2,3$, and $\bar t_1 \bar t_2 \bar t_3 \bar t_4$.  Hence, by (\ref{cohomGENS}), the generators of $\H^*(\B_{\G})$ can be chosen as $\bar y_i \ceq \sigma_i (\bar t^2) \in \H^{4i}$, $i=1,2,3$, and $\bar y_4 \ceq \bar t_1 \bar t_2 \bar t_3 \bar t_4 \in \H^8$.

$W_{\UU_2}$ is the dihedral group of order twelve.  It acts trivially on $u$, while the action on the root system $\G_2$ is by rotations of $\frac{\pi}{3}$ and by reflections through the horizontal axis.  Therefore, given our description of the root system of $\G_2$ above, $W_{\UU_2}$ acts on $\H^*(\B_{\T_\UU})$ via permutations in $\bar s_i$ and a simultaneous sign change of all $\bar s_i$.  On the other hand, $W_{\UU_3} = \Z_2 \x W_{\UU_2}$, where the $\Z_2$ factor acts trivially on the $s_i$ and by a sign change on $u$.  Thus elements of $\H^*(\B_{\T_\UU})$ which are invariant under $W_{\UU_k}$ are given by sums and products of $\bar u^{k-1}$ and the elementary symmetric polynomials $\sigma_2 (\bar s) \ceq \sigma_2 (\bar s_1, \bar s_2, \bar s_3)$ and $\sigma_i (\bar s^2) \ceq \sigma_i (\bar s_1^2, \bar s_2^2, \bar s_3^2)$, $i=1,2,3$.  However, since $\sum s_i = 0$, Lemma \ref{nicetrick} shows that a basis for $\H^*(\B_{\T_\UU})^{W_{\UU_k}}$ is given by $\bar u^{k-1}$ and the symmetric polynomials $\sigma_2 (\bar s)$ and $\sigma_3 (\bar s^2)$.  By (\ref{cohomGENS}), generators of $\H^*(\B_{\UU_k})$ are given by $\bar w_k \ceq \bar u^{k-1} \in \H^{2(k-1)}$, $\bar x_1 \ceq \sigma_2 (\bar s) \in \H^4$ and $\bar x_2 \ceq \sigma_3 (\bar s^2) \in \H^{12}$.

Consider $\iota_k \co \UU_k \hra \GG$ as above.  By (\ref{maxtorus}), $\iota_k|_{\T_\UU}\co \T_\UU \lra \T_{\GG}$ has (in coordinates) the form $\iota_k|_{\T_\UU} (u, s_1,s_2,s_3) = ((0, 0, 0, u), (0,s_1,s_2,-s_3))$.  Hence
        \begin{align*}
            (\iota_k|_{\T_\UU})^* \co \H^1(\T_{\GG}) &\lra \H^1(\T_\UU)\\
                            t_i \ox 1 &\lmt 0, \ \ i = 1, 2, 3 \\
                            t_4 \ox 1 &\lmt u \\
                            1 \ox t_1 &\lmt 0\\
                            1 \ox t_i &\lmt s_{i-1}, \ \ i = 2,3\\
                            1 \ox t_4 &\lmt -s_3
        \end{align*}
from which it follows:
        \begin{align*}
            (\B_{\iota_k|_{\T_\UU}})^* \co \H^2(\B_{\T_{\GG}}) &\lra \H^2(\B_{\T_\UU})\\
                            \bar t_i \ox 1 &\lmt 0, \ \ i = 1, 2, 3 \\
                            \bar t_4 \ox 1 &\lmt \bar u \\
                            1 \ox \bar t_1 &\lmt 0\\
                            1 \ox \bar t_i &\lmt \bar s_{i-1}, \ \ i = 2,3\\
                            1 \ox \bar t_4 &\lmt - \bar s_3
        \end{align*}
Therefore, by (\ref{cohomdiagram}) we have:
        \begin{align}
            \nonumber   (\B_{\iota_k})^* \co \H^*(\B_{\GG}) &\lra \H^*(\B_{\UU_k})\\
            \nonumber       \bar y_1 \ox 1 &\lmt \bar u^2 \\
            \label{Bg}      \bar y_i \ox 1 &\lmt 0, \ \ i = 2,3,4\\
            \nonumber       1 \ox \bar y_i &\lmt \sigma_i (\bar s^2) \ \ i = 1, 2,3\\
            \nonumber       1 \ox \bar y_4 &\lmt 0
        \end{align}
In particular, $(\B_{\iota_k})^* (1 \ox \bar y_1) = - 2 \bar x_1$, $(\B_{\iota_k})^* (1 \ox \bar y_2) = \bar x_1^2$ (by Lemma \ref{nicetrick}) and $(\B_{\iota_k})^* (1 \ox \bar y_3) = \bar x_2$.
%\ms

We are now in a position to compute the first Pontrjagin class of $\G \bq \UU_k$.  Recall that the positive roots of $\G = \SO(8)$ are $t_i \pm t_j$, $1 \leq i < j \leq 4$.  Hence:
        \begin{align*}
        \sum_{\alpha_i \in \Delta_{\G}^+} \bar \alpha_i^2 &= \sum_{1 \leq i < j \leq 4} \left((\bar t_i - \bar t_j)^2 + (\bar t_i + \bar t_j)^2\right) \\
                        &= 2 \sum_{1 \leq i < j \leq 4} (\bar t_i^2 + \bar t_j^2)\\
                        &= 6 \sum_{i=1}^4 \bar t_i^2 \\
                        &= 6 \bar y_1
        \end{align*}
Since $\SO(2)$ has no (positive) roots and $\SO(3)$ has only one positive root, namely $u$, we may denote the positive root of $\SO(k)$ by $(k-2)u$, for $k = 2,3$, respectively.  From our earlier description of the roots of $\G_2$, the positive roots of $\UU_k = \SO(k) \x \G_2$ are
        $$(k-2)u,\ s_1,\ s_2,\ - s_3,\ s_1 - s_3,\ s_2 - s_1,\ s_2 - s_3$$
where $\sum s_i = 0$.  Then:
        \begin{align*}
        \sum_{\beta_j \in \Delta_{\UU_k}^+} \bar \beta_j^2 &= (k-2)\bar u^2 + \bar s_1^2 + \bar s_2^2 + \bar s_3^2 + (\bar s_1 - \bar s_3)^2 + (\bar s_2 - \bar s_1)^2 + (\bar s_2 - \bar s_3)^2\\
                        %&= (k-2)\bar u^2 + \sigma_1 (\bar s^2) \!+\! (\bar s_1^2 + \bar s_3^2 - 2\bar s_1 \bar s_3) \!+\! (\bar s_1^2 + \bar s_2^2 - 2\bar s_1 \bar s_2) \!+\! (\bar s_2^2 + \bar s_3^2 - 2\bar s_2 \bar s_3)\\
                         &= (k-2)\bar u^2 + 3 \sigma_1 (\bar s^2) - 2 \sigma_2 (\bar s)\\
                         &= (k-2)\bar u^2  -8 \sigma_2 (\bar s) \qquad \qquad \ \textrm{by Lemma \ref{nicetrick}}\\
                         &= (k-2)\bar u^2 -8 \bar x_1
        \end{align*}
Note that, by (\ref{equality}) and  (\ref{Bg}),
\begin{align*}
\varphi_{\G}^* (\bar y_1) &=
\vphi_{\UU_k}^*((\B_{\iota_k})^*(\bar y_1 \ox 1)) =
\vphi_{\UU_k}^*(\bar u^2) \ \ \ \ \ {\rm and}\\
- 2 \varphi_{\UU_k}^* (\bar x_1) &=
\vphi_{\UU_k}^*((\B_{\iota_k})^*(1 \ox \bar y_1)) =
\vphi_{\UU_k}^*((\B_{\iota_k})^*(\bar y_1 \ox 1)) =
\vphi_{\UU_k}^*(\bar u^2).
\end{align*}
It follows now from (\ref{firstPont}) that:
\begin{align*}
            \pp_1 (\G \bq \UU_k) &= \pp_1(\alpha_{\G}) - \pp_1(\alpha_{\UU_k})\\
                &= 6 \vphi_{\G}^* (\bar y_1) - (k-2) \vphi_{\UU_k}^*(\bar u^2) + 8 \vphi_{\UU_k}^*(\bar x_1) \\
                &= 6 \vphi_{\UU_k}^*(\bar u^2) - (k-2) \vphi_{\UU_k}^*(\bar u^2) - 4 \vphi_{\UU_k}^*(\bar u^2) \\
                &= (4-k) \vphi_{\UU_k}^*(\bar u^2) \in \H^4(\G \bq \UU_k)
        \end{align*}
It remains only to show that $\vphi_{\UU_k}^*(\bar u^2)$ is a generator of $\H^4(\G \bq \UU_k)$.  To achieve this, consider the Serre spectral sequence for the fibration $\G \lra \G \bq \UU_k \lra \B_{\UU_k}$.

Suppose $k = 2$.  Notice that $\< \bar w_2 \> = \H^2(\B_{\UU_2}) = \E^{2,0}_2$ will survive until $\E_{\infty}$ since $\H^*(\G)$ contains no elements of degree $1$.

On the other hand, suppose now that $k = 3$.  By (\ref{Bg}) the generator $y_1 \in \H^3 (\G) = \E^{0,3}_2 = \E^{0,3}_4$ gets mapped under $d_4$ to
$$
d_4(y_1) = (\B_{\iota_3})^*(\bar y_1 \ox 1) - (\B_{\iota_3})^*(1 \ox \bar y_1) = \bar u^2 + 2 \bar x_1 \in \E^{4,0}_4 = \E^{4,0}_2 = \H^4(\B_{\UU_3}).
$$
The generators $\bar w_3 = \bar u^2$ and $\bar x_1$ of $\H^4(\B_{\UU_3})$ are both mapped to zero by $d_4$.  Thus, the $\E^{4,0}_5$ term is a $\Z_p$ generated by $\bar w_3 = \bar u^2$ and survives to $\E_\infty$.

Recall that the classifying map $\vphi_{\UU_k}^*$ is the edge homomorphism
        $$\vphi_{\UU_k}^* \co \H^i (\B_{\UU_k}) = \E^{i,0}_2 \twoheadrightarrow \E^{i,0}_\infty \hra \H^i (\G \bq \UU_k).$$
By the discussion above $\vphi^*_{\UU_k} (\bar w_k) = \vphi^*_{\UU_k} (\bar u^{k-1}) \neq 0$, $k = 2,3$.  Hence $\vphi^*_{\UU_k} (\bar w_k)$ is a generator of $\H^{2(k-1)}(\G \bq \UU_k) = \Z_p$.  If $k = 3$ we are done.  When $k = 2$, we know from Theorem \ref{cohom} that $\H^4(\M^{13})$ is generated by the square of a generator of $\H^2(\M^{13})$.  That is, $\vphi^*_{\UU_2} (\bar w_2^2) = \vphi^*_{\UU_2} (\bar u^2)$ is a generator of $\H^4(\M^{13})$.
\end{proof}

\begin{rmk}
Since $\H^8$ and $\H^{12}$ are trivial for each of the manifolds $\M^{13}$ and $\N^{11}$, we have in fact computed their total Pontrjagin classes $\pp = 1 + \pp_1$ in $\Z_p$ coefficients.
\end{rmk}

\begin{rmk}
\label{Pontrmk}
In terms of integral cohomology, the proposition tells us only that $\pp_1(\M^{13})$ and $\pp_1(\N^{11})$ are non-trivial and not divisible by any primes $p \geq 3$. Thus $\pp_1(\M^{13})$ and $\pp_1(\N^{11})$ have the form $\pm 2^\ell \in \Z = \H^4$, for some $\ell \in \NN \cup \{0\}$.
\end{rmk}

\begin{cor}
$\N^{11}$ is not homeomorphic to $\sph^4 \x \sph^7$.
\end{cor}
\begin{proof}
By the previous remark, in integral coefficients $\pp_1 (\N^{11}) = \pm 2^\ell \in \Z = \H^4(\N^{11}; \Z)$, for some $\ell \in \NN \cup \{0\}$.  However, since all (integral) Pontrjagin classes for spheres are trivial and integral Pontrjagin classes are homeomorphism invariants, $\N^{11}$ cannot be homeomorphic to a product of spheres.
\end{proof}

Since $\pp_1(\CP^3 \x \sph^7) = 4\gamma^2$, where $\gamma$ generates $\H^2(\CP^3 \x \sph^7; \Z)$, we are unable to distinguish $\M^{13}$ and $\CP^3 \x \sph^7$ using the proposition.  We need to explicitly compute the integral Pontrjagin class of $\M^{13}$.  We can accomplish this by ``hot-wiring'' the technique for computing Pontrjagin classes in the absence of torsion in the cohomology groups.

Before we begin we establish two topological statements which will be used in the proof of Theorem \ref{hotwire}.  From now on we will always assume that our cohomology groups have integral coefficients, and by spectral sequence we will always mean Serre spectral sequence.

\begin{prop}
\label{BG2}
$\B_{\G_2}$, the classifying space of $\G_2$, has low dimensional integral cohomology groups $\H^1 = \H^2 = \H^3 = \H^5 = 0$ and $\H^4 = \Z$ with generator $\bar x =  \sigma_2(\bar s)$, where $\sigma_2(\bar s) \ceq \sigma_2(\bar s_1, \bar s_2, \bar s_3)$, $\sum \bar s_i = 0$, and $\bar s_i \in \H^2 (\B_{\T_{\G_2}})$, $i=1,2,3$, are the transgressions of the elements $s_i \in \H^1 (\T_{\G_2})$, $i=1,2,3$, which span the integral lattice of $\G_2$.
\end{prop}

\begin{proof}
Consider the universal bundle $\G_2 \lra \E_{\G_2} \lra \B_{\G_2}$ where $\E_{\G_2}$ is contractible.  From Whitehead \cite[Theorem 5.17]{Wh}, we know that $\H^j (\G_2) = 0$, $j=1,2,4,5$, and $\H^3 (\G_2) = \Z$.  Let $x$ be a generator of $\H^3 (\G_2)$.  Since $\E_{\G_2}$ is contractible all entries in the spectral sequence for the fibration $\G_2 \lra \E_{\G_2} \lra \B_{\G_2}$ must get killed off.  Since $d_4 \co \E_4^{0,3} \lra \E_4^{4,0}$ is the only possible non-trivial differential with domain $\E_4^{0,3}$ it must map $x \in \H^3 (\G_2)$ to a generator $\bar x$ of $\H^4 (\B_{\G_2})$, and so $\H^4 (\B_{\G_2}) = \Z$.  Similarly it is clear from the spectral sequence that $\H^j (\B_{\G_2}) = 0$ for $j = 1,2,3,5$.

Now consider the fibration $\sph^6 = \G_2 / \SU(3) \lra \B_{\SU(3)} \lra \B_{\G_2}$.  The spectral sequence associated to this fibration shows that $\bar x \in \E_2^{4,0} = \H^4 (\B_{\G_2})$ survives to $\E_\infty$.  Thus, since there are no other non-zero entries on the corresponding diagonal in $\E_\infty$, we see that $\H^4 (\B_{\G_2}) = \H^4 (\B_{\SU(3)})$.  Recall that $\H^* (\B_{\SU(3)})$ is a polynomial algebra generated by the elementary symmetric polynomials $\sigma_i (\bar s) = \sigma_i (\bar s_1, \bar s_2, \bar s_3)$, $i=2,3$, in the transgressions $\bar s_j$ of $s_j \in \H^1 (\T_{\SU(3)})$, $j = 1,2,3$, where the $s_j$ span the integral lattice of $\SU(3)$.  Note that $\sum s_j = 0$, $\T_{\G_2} = \T_{\SU(3)}$ and $\deg(\sigma_i (\bar s)) = 2i$.  Therefore $\H^4(\B_{\G_2})$ is generated by $\sigma_2 (\bar s)$ as desired.
\end{proof}

\begin{prop}
\label{SO(8)/G2}
The low dimensional integral cohomology groups of the manifold $\SO(8)/\G_2 = (\sph^7 \x \sph^7)/\Z_2$ are $\H^j (\SO(8)/\G_2) = \H^j (\rp^7)$, $0 \leq j \leq 6$.
\end{prop}

\begin{proof}
Consider the spectral sequence for the fibration
        $$\rp^7 = \SO(7)/\G_2 \lra \SO(8)/\G_2 \lra \SO(8)/\SO(7) = \sph^7.$$
Recall:
        $$\H^j (\rp^7) = \begin{cases}
                        \Z & \text{if $j = 0,7$}\\
                        \Z_2 & \text{if $j = 2,4,6$}\\
                        0 & \text{if $j = 1,3,5$}
                        \end{cases}$$
It is clear that each $\E_2^{0,j} = \H^j (\rp^7)$, $j \leq 5$, survives to $\E_\infty$.  For $\E_2^{0,6} = \H^6 (\rp^7) = \Z_2$ notice that there are no non-trivial homomorphisms $\Z_2 \lra \Z$ and so the differential $d_7 \co \E_7^{0,6} = \Z_2 \lra \E_7^{7,0} = \Z$ must be trivial.  Therefore $\E_2^{0,6} = \H^6 (\rp^7)$ also survives to $\E_\infty$.  Since there are no other non-zero entries on the corresponding diagonals we get the desired result.
\end{proof}

We are now ready to complete the proof of Theorem \ref{thmA}\ref{QPnew1}.

\begin{thm}
\label{hotwire}
The first integral Pontrjagin class of $\M^{13} = \SO(8) \bq (\sph^1 \x \G_2)$ is given by
        $$|\pp_1(\M^{13})| = 8 z^2$$
where $z$ is a generator of $\H^2 (\M^{13}; \Z) = \Z$.  In particular, $\M^{13}$ is not homeomorphic to $\CP^3 \x \sph^7$.
\end{thm}

\begin{proof}
Consider once again diagram (\ref{diagram}), with $\G = \SO(8)$ and $\UU = \sph^1 \x \G_2$ such that $\G \bq \UU = \M^{13}$.  In the proof of Proposition \ref{Pontclass13} we followed the usual techniques of Borel and Hirzebruch \cite{BH}, Eschenburg \cite{E2} and Singhof \cite{Si} when there is no torsion in cohomology, namely we computed $\B_\iota$ and $\B_\Delta$ and then used the fact that the diagram commutes in order to compute the $\Z_p$ Pontrjagin class, for odd primes $p$.  However, since $\SO(8)$ and $\G_2$ have torsion in integral cohomology, we need to adopt a different approach in order to compute the integral Pontrjagin class.  Since $\H^8 (\M^{13}) = \H^{12} (\M^{13}) = 0$ we can restrict our attention to the first integral Pontrjagin class $\pp_1 (\M^{13}) \in \H^4 (\M^{13})$.  The key idea to be taken from the proof of Proposition \ref{Pontclass13} is that we computed the first Pontrjagin classes of some vector bundles over $\B_{\G}$ and $\B_\UU$, then pulled them back to $\M^{13}$ under the classifying maps $\vphi_{\G}$ and $\vphi_\UU$ respectively.  As it turns out, the first Pontrjagin classes of these vector bundles over $\B_{\G}$ and $\B_\UU$ are the same in integral coefficients as in $\Z_p$ coefficients $p \geq 3$.  Our strategy, therefore, is to compute the maps $\vphi_\UU^* \co \H^4 (\B_\UU) \lra \H^4 (\M^{13})$ and $\vphi_{\G}^* \co \H^4 (\B_{\G}) \lra \H^4 (\M^{13})$ and pull back the respective first Pontrjagin classes.

As a first step in computing $\vphi_\UU^* \co \H^4 (\B_\UU) \lra \H^4 (\M^{13})$ we notice that $\H^* (\UU) = \H^*(\sph^1) \ox \H^*(\G_2)$ and $\H^*(\B_\UU) = \H^*(\B_{\sph^1}) \ox \H^*(\B_{\G_2})$ since $\H^*(\sph^1)$ and $\H^*(\B_{\sph^1})$ are torsion-free.  Therefore
        $$\H^j (\UU) = \begin{cases}
                    \Z & \text{if $j = 0$}\\
                    \Z = \<w\> & \text{if $j = 1$}\\
                    \Z = \<x\> & \text{if $j = 3$}\\
                    0 & \text{if $j = 2,4,5$}
                    \end{cases}$$
where $w$ is a generator of $\H^1 (\sph^1)$ and $x$ is a generator of $\H^3 (\G_2)$,  and applying Proposition \ref{BG2}
        $$\H^j (\B_\UU) = \begin{cases}
                    \Z & \text{if $j = 0$}\\
                    \Z = \<\bar w\> & \text{if $j = 2$}\\
                    \Z \oplus \Z = \<\bar w^2\> \oplus \<\bar x\>& \text{if $j = 4$}\\
                    0 & \text{if $j = 1,3,5$}
                    \end{cases}$$
where $\bar w$ is the transgression of $w$ resulting from the spectral sequence for the universal bundle of $\sph^1$ and generates $\H^2 (\B_{\sph^1})$ (hence generates $\H^* (\B_{\sph^1}) = \Z[\bar w]$), and $\bar x$ is the transgression of $x$ resulting from the spectral sequence for the universal bundle of $\G_2$ and generates $\H^4 (\B_{\G_2})$.

Recall that $\vphi_\UU \co \G \bq \UU \lra \B_\UU$ is the classifying map since we have the following diagram of principal $\UU$-bundles
        $$\xymatrix{ \UU \ar[r] \ar[d] & \UU \ar[d] \\
                        \G \x \E_\UU \ar[r]^(.6){\pi_2} \ar[d] & \E_\UU \ar[d] \\
                        \G \x_\UU \E_\UU \ar[r]_(.6){\pi_2} & \B_\UU
        }$$
where $\pi_2$ denotes projection onto the second factor and $\UU \lra \E_\UU \lra \B_\UU$ is the universal bundle.  Since $\E_\UU$ is contractible, projection onto the first factor gives homotopy equivalences $\G \x \E_\UU \simeq \G$ and $\G \x_\UU \E_\UU \simeq \G \bq \UU$.  The resulting map $\G \bq \UU \lra \B_\UU$ is $\vphi_\UU$.  Hence $\vphi_\UU$ is the classifying map.  Therefore, up to homotopy, we may consider the following commutative diagram of fibrations
        $$\xymatrix{ \UU \ar[r] \ar[d] & \UU \ar[d] \\
                        \G \ar[r] \ar[d] & \E_\UU \ar[d] \\
                        \G \bq \UU \ar[r]_{\vphi_\UU} & \B_\UU
        }$$
Consider first the spectral sequence for the fibration on the left.  Recall that $\H^*(\M^{13}) = \H^*(\CP^3 \x \sph^7)$.  Hence:
        $$\H^j (\G \bq \UU) = \begin{cases}
                    \Z & \text{if $j = 0$}\\
                    \Z = \<z\> & \text{if $j = 2$}\\
                    \Z = \<z^2\> & \text{if $j = 4$}\\
                    0 & \text{if $j = 1,3,5$}
                    \end{cases}$$
Since $\G = \SO(8)$ we have from \v Cadek, Mimura and Van\v zura \cite{CMV}:
        $$\H^j (\G) = \begin{cases}
                    \Z & \text{if $j = 0$}\\
                    0 & \text{if $j = 1$}\\
                    \Z_2 = \<r\> & \text{if $j = 2$}\\
                    \Z = \<y\> & \text{if $j = 3$}\\
                    \Z_2 = \<r^2\> & \text{if $j = 4$}
                    \end{cases}$$
Since $\H^1 (\G) = 0$ we see that $d_2 \co \E_2^{0,1} = \<w\> \lra \E_2^{2,0} = \<z\>$ must have trivial kernel, i.e. $d_2 (w) = k z$ for some $k \in \Z$, $k \neq 0$.  Then $\E_3^{0,2} = \<z\>/\<k z\>$ survives to $\E_\infty$ and since $\H^2 (\G) = \Z_2$ we must therefore have $k = \pm 2$, i.e. $d_2 (w) = \pm 2z$.

On the other hand, the spectral sequence shows that on the $\E_4$-page we have the differential $d_4 \co \E_4^{0,3} = \<x\> \lra \E_4^{0,4} = \<z^2\> / \<2 z^2\>$.  However, since $\H^3 (\G) = \Z$ and $\H^4 (\G) = \Z_2$, we must have $d_4 (x) = 0 \in \<z^2\> / \<2 z^2\>$.

Since $\E_\UU$ is contractible it is clear from the spectral sequence for the fibration on the right that $d_2 \co \E_2^{0,1} = \< w \> \lra \E_2^{2,0} = \<\bar w\>$ is an isomorphism with $d_2(w) = \bar w$, and $d_4 \co \E_4^{0,3} = \<x\> \lra \E_4^{4,0} = \<\bar w^2\> \oplus \<\bar x\>$ is given by $d_4 (x) = \bar x$.

By naturality of the spectral sequence we thus have for the left-hand fibration that $d_2 (w) = \vphi_\UU^* (\bar w) \in \<z\>$ and $d_4(x) = \vphi_\UU^* (\bar x) \in \<z^2\> / \<2 z^2\>$.  Therefore, since we have already shown that $d_2 (w) = \pm 2z \in \<z\>$ and $d_4 (x) = 0 \in \<z^2\> / \<2 z^2\>$, we find
        \begin{align}
        \nonumber
            \vphi_\UU^* (\bar w) &= \pm 2z \in \H^2 (\G \bq \UU) = \<z\> \ \ \ {\rm and}\\
        \label{phiUx}
            \vphi_\UU^* (\bar x) &= 2kz^2 \in \H^4 (\G \bq \UU) = \<z^2\>, \ \ \ \textrm{for some $k \in \Z$}.
        \end{align}
We now turn our attention to computing $\vphi_{\G}^* \co \H^4 (\B_{\G}) \lra \H^4 (\M^{13})$.  In order to show that $\vphi_{\G} \co \G \bq \UU \lra \B_{\G}$ is the classifying map consider the commutative diagram of principal $\G$-bundles
        $$\xymatrix{  \G \ar[r] \ar[d] & \G \ar[d] \\
                        \GG \x_\UU \E_{\GG} \ar[r] \ar[d] & \GG \x_{\GG} \E_{\GG} \ar[d]\\
                       (\Delta \G \backslash \GG) \x_\UU \E_{\GG} \ar[r]_{\vphi_{\G}} & (\Delta \G \backslash \GG) \x_{\GG} \E_{\GG}
        }$$
Since $\GG \x_{\GG} \E_{\GG} = \E_{\GG}$ and $(\Delta \G \backslash \GG) \x_{\GG} \E_{\GG} = \G \x_{\GG} \E_{\GG} = \B_{\Delta \G}$ we see that the fibration on the right-hand side is the universal bundle for $\G$.  On the left-hand side we have $(\Delta \G \backslash \GG) \x_\UU \E_{\GG} = \G \x_\UU \E_{\GG}$, and projection onto the first factor gives homotopy equivalences $\GG \x_\UU \E_{\GG} \simeq \GG/\UU$ and $\G \x_\UU \E_{\GG} \simeq \G \bq \UU$.  Thus up to homotopy the diagram becomes
        $$\xymatrix{  \G \ar[r] \ar[d] & \G \ar[d] \\
                        \GG/\UU \ar[r] \ar[d] & \E_{\GG} \ar[d]\\
                       \G \bq \UU \ar[r]_{\vphi_{\G}} & \B_{\G}
        }$$
as desired.  Recall that $\H^3 (\G) = \<y\>$.  The cohomology of $\B_{\G}$ is described in Brown \cite{Br} and Feshbach \cite{F}, but for our purposes we need only that:
        $$\H^j (\B_{\G}) = \begin{cases}
                    \Z & \text{if $j = 0$}\\
                    0 & \text{if $j = 1,2$}\\
                    \Z_2 & \text{if $j = 3$}\\
                    \Z & \text{if $j = 4$}\\
                    \Z_2 & \text{if $j = 5$}
                    \end{cases}$$
Whilst proving Proposition 3.6 in \cite{GZ2} Grove and Ziller showed that, since $\E = \E_{\GG}$ is contractible, in the spectral sequence for the bundle $\G \lra \E \lra \B_{\G}$ the differential $d_4 \co \E_4^{0,3} = \<2y\> \lra \E_4^{4,0} = \H^4(\B_{\G})$ is an isomorphism, i.e. $2y$ gets mapped to a generator $\bar y$ of $\H^4 (\B_{\G}) = \Z$.  This follows from the facts that $\E_2^{2,2} = \Z_2$ (by the Universal Coefficient Theorem) and that $d_2 \co \E_2^{0,3} = \<y\> \lra \E_2^{2,2} = \Z_2$ must be onto.

Therefore naturality of the spectral sequence implies that $d_4 (2y) = \vphi_{\G}^* (\bar y)$ in the spectral sequence for the left-hand fibration $\G \lra \GG/\UU \lra \G \bq \UU$, where $\H^3 (\G) = \<y\>$ and $\H^4(\B_{\G}) = \<\bar y\>$.

In order to determine the exact value of $\vphi_{\G}^*(\bar y) \in \H^4(\G \bq \UU)$ we need to examine the spectral sequence for the left-hand fibration.  First we must compute the cohomology of $\GG/\UU$ in low-dimensions.  Recall that $\GG/\UU = V_{8,6} \x \SO(8)/\G_2$, where $V_{8,6}$ is the Stiefel manifold $\SO(8)/\SO(2)$.  From \v Cadek, Mimura and Van\v zura \cite{CMV} we find:
        $$\H^j (V_{8,6}) = \begin{cases}
                    \Z & \text{if $j = 0,2$}\\
                    0 & \text{if $j = 1,3,5$}\\
                    \Z_2 & \text{if $j = 4$}
                    \end{cases}$$
In Proposition \ref{SO(8)/G2} we computed the low dimensional cohomology groups of $\SO(8)/\G_2$.  From the general K\"unneth formula for cohomology (Spanier \cite[Theorem 11, Section 5]{Sp}) it follows:
        $$\H^j (\GG/\UU) = \begin{cases}
                    \Z & \text{if $j = 0$}\\
                    0 & \text{if $j = 1,3$}\\
                    \Z \oplus \Z_2 & \text{if $j = 2$}\\
                    \Z_2^3 & \text{if $j = 4$}\\
                    \Z_2 & \text{if $j = 5$}
                    \end{cases}$$
Since $\H^4 (\GG/\UU) = \Z_2^3$, in the spectral sequence for $\G \lra \GG/\UU \lra \G \bq \UU$ the differential $d_2 \co \E_2^{0,3} = \H^3(\G) = \<y\> \lra \E_2^{2,2} = \Z_2$ must be trivial, i.e. $\E_2^{2,2} = \Z_2$ must survive to $\E_\infty$.  It thus follows that $\E_2^{i,j} = \E_3^{i,j} = \E_4^{i,j}$ for $i \leq 5, j \leq 4$.  Since $\H^3(\GG/\UU) = 0$ the differential $d_4 \co \E_4^{0,3} = \<y\> \lra \E_4^{4,0} = \H^4(\G \bq \UU) = \<z^2\>$ must be given by $d_4 (y) = n z^2$ for some non-zero $n \in \Z$.  On the other hand, since $\H^4 (\GG/\UU) = \Z_2^3$, $\E_4^{0,4} = \E_4^{2,2} = \Z_2$ and $\E_4^{1,3} = \E_4^{3,1} = 0$, the filtration for the spectral sequence shows that $n = \pm 2$, i.e. $d_4 (y) = \pm 2 z^2$. But we have already shown that $d_4 (2y) = \vphi_{\G}^* (\bar y)$.  Therefore
        $$\vphi_{\G}^*(\bar y) = \pm 4z^2 \in \H^4 (\G \bq \UU) = \<z^2\>.$$
Furthermore, while proving Lemma 5.4 in \cite{GZ2} Grove and Ziller showed that, by considering the spectral sequences of the fibrations $\SO(8)/\SO(3) \lra \B_{\SO(3)} \lra \B_{\SO(8)}$ and $\SO(3)/\SO(2) \lra \B_{\SO(2)} \lra \B_{\SO(3)}$, we can let $\bar y = \sigma_1 (\bar t^2) = \sigma_1 (\bar t_1^2, \bar t_2^2, \bar t_3^2, \bar t_4^2)$, where $(t_1, \dots, t_4)$ are the coordinates of a maximal torus $\T_{\G}$ of $\G$ and by abuse of notation we identify $t_i \in \H^1 (\T_{\G})$ with $\bar t_i \in \H^2 (\B_{\T_{\G}})$ via transgression.

We are now in a position to compute the first Pontrjagin class of $\M^{13} = \G \bq \UU$.  We have already shown that, despite having torsion in cohomology, $\H^4(\B_{\G}) \cong \H^4(\B_{\T_{\G}})^{W_{\G}}$ and $\H^4(\B_\UU) \cong \H^4(\B_{\T_\UU})^{W_\UU}$ since the generators are $\bar y = \sigma_1(\bar t^2)$ and $\bar x = \sigma_2(\bar s)$ respectively.  Therefore, in the present setting, equation (\ref{firstPont}) for $\pp_1$ is still valid for integral coefficients.

As in the proof of Proposition \ref{Pontclass13}, $\pp_1 (\alpha_{\G}) = 6 \vphi_{\G}^*(\bar y)$.   But $\varphi_{\G}^* (\bar y) = \pm 4 z^2$.  Hence $\pp_1 (\alpha_{\G}) = \pm 24 z^2 \in \H^4 (\G \bq \UU)$.

Similarly, from the proof of Proposition \ref{Pontclass13} we have $\pp_1 (\alpha_\UU) = - 8 \vphi_\UU^*(\bar x)$.  Thus, since $\varphi_\UU^* (\bar x) = 2k z^2$ by (\ref{phiUx}), $\pp_1 (\alpha_\UU) = - 16 k z^2 \in \H^4 (\G \bq \UU)$, for some $k \in \Z$.

Therefore, by (\ref{firstPont}),
        \begin{align*}
            \pp_1 (\tau) &= \pp_1(\alpha_{\G}) - \pp_1(\alpha_\UU)\\
                &= 8(2k \pm 3)z^2 \in \H^4(\G \bq \UU).
        \end{align*}
By Proposition \ref{Pontclass13} we know that $\pp_1 (\tau) = \pp_1 (\G \bq \UU)$ is divisible only by $2$.  Therefore we must have $2k \pm 3 = \pm 1$, which implies $\pp_1(\G \bq \UU) = \pm 8 z^2$ as desired.
\end{proof}

%%%%%%%%%%%%%%%%%%%%   End of main body of article
%
%                             References
%
%   BiBTeX users uncomment the following line:
%
%\bibliographystyle{gtart}
%

\end{document}